\documentclass[11pt,a4paper]{amsart}

\newcommand{\figdir}{}

\usepackage[foot]{amsaddr}

\usepackage[margin=3cm]{geometry}
\usepackage{amsmath,amsfonts,amssymb,amsthm}
\usepackage{graphicx}
\usepackage{clrscode3e}
\usepackage{enumerate}

\newcommand{\pdpd}[2]{\frac{\partial #1}{\partial #2}}

\newtheorem{theorem}{Theorem}

\newcommand{\R}{{\mathbb{R}}}

\newcommand{\NIL}{\const{nil}}

\title[An improved sweeping domain decomposition preconditioner]%
{An improved sweeping domain decomposition preconditioner for the
  Helmholtz equation}
\author{Christiaan C. Stolk}
\email{C.C.Stolk@uva.nl}
\address{University of Amsterdam, 
Korteweg-de Vries Institute for Mathematics, 
P.O.Box 94248, 1090 GE Amsterdam, The Netherlands}

\begin{document}
\begin{abstract}
\noindent
In this paper we generalize and improve a recently developed domain
decomposition preconditioner for the iterative solution of discretized
Helmholtz equations.  We introduce an improved method for transmission
at the internal boundaries using perfectly matched layers.
Simultaneous forward and backward sweeps are introduced, thereby
improving the possibilities for parallellization. Finally, the method
is combined with an outer two-grid iteration. The method is studied
theoretically and with numerical examples. It is shown that the
modifications lead to substantial decreases in computation time and
memory use, so that computation times become comparable to that of the
fastests methods currently in the literature for problems with up to
$10^8$ degrees of freedom.
\end{abstract}

\maketitle


\section{Introduction}

The linear systems resulting from discretizing the high-frequency
Helmholtz equation have been a challenge for mathematicians for a long
time \cite{Erlangga2008,ErnstGander2011}.  A class of methods that
recently gained much attention is that of sweeping domain
decomposition preconditioners and related methods
\cite{EngquistYing2011,Stolk2013,PoulsonEtAl2013,ChenXiang2013,
  VionGeuzaine2014,ZepedaNunezDemanet2016_Polarized_traces}.  In this paper we consider
the improvement and generalization of one such method, namely a double
sweep method using the perfectly matched layer (PML) at the
interfaces, described in \cite{Stolk2013}.

To be specific, we consider the Helmholtz equation in two and three
dimensions. In two dimensions it reads
\begin{equation}
  -\partial_{xx}^2 u(x,y) - \partial_{yy}^2u(x,y) - k(x,y)^2 u(x,y) 
  = f(x,y) ,
\end{equation}
where $k(x,y) = \frac{\omega}{c(x,y)}$, with $c(x,y)$ the wave speed.
The computational domain is assumed to be a rectangle that is
truncated using perfectly matched layers, or classical damping layers.
We consider finite difference or finite element discretizations on
regular meshes that result in a compact stencil, i.e.\ a $3 \times 3$
or $3 \times 3 \times 3$ square or cubic stencil depending on the
dimension. Accurate discretizations of this type are possible, see
e.g.\ \cite{TurkelEtAl2013,Stolk2016_Dispersion_minimizing_scheme}.  Thus we
generalize the results of \cite{Stolk2013} involving second order
finite differences.

Domain decomposition methods for the Helmholtz equation typically
follow, to an extent depending on numerical approximations, the
principles that:
\begin{equation} \label{eq:conditions}
\begin{minipage}{13cm}
\begin{enumerate}[(i)]
\item
the boundary conditions at the subdomain interfaces should be non-reflecting; 
\item
if $\Omega^{(j-1)}$ and $\Omega^{(j)}$ are neighboring subdomains then
the outgoing wave field from $\Omega^{(j-1)}$ should equal the  
incoming wave field in $\Omega^{(j)}$ at the joint boundary and vice versa.
\end{enumerate}
\end{minipage}
\end{equation}
The use of Robin or numerical absorbing boundary conditions at the
interfaces is one way to do this, see e.g.\
\cite{GanderHalpernMagoules2007,VionGeuzaine2014} and
references. Another way is using PML boundary layers
\cite{Stolk2013,SchadleEtAl2007}, the method we will use here (in
modified form).

{\em Double sweep domain decomposition} is distinguished from other
domain decomposition methods by the ordering of the subdomain
solves. Here the subdomains are chosen as parallel slices of the
original domain, say numbered from 1 to $J$. The subdomain solutions
are computed first for $j=1,\ldots,J$ subsequently, this is called the
forward sweep, and then for $j=J,J-1,\ldots, 1$ subsequently, called
the backward sweep. In this way information can propagate over the
entire domain in one preconditioner application. The condition that
information can propagate over at least an $O(1)$ part of the domain
is necessary to achieve a good approximation of the true solution.

In this paper we consider three modifications to the method of
\cite{Stolk2013}. The first concerns the transmission of information
between neighboring subdomains using PML layers. It was observed in
\cite{Stolk2013} that at the onset of the PML layer the field is
approximately outgoing and that in the next subdomain, a similar
ingoing field can be reproduced using a planar source proportional to
the outgoing field.  We modify the way this is done compared to
\cite{Stolk2013}. The new method is more generally applicable and
prevents the planar source radiating into the added absorbing layer,
which is an advantage because these layers are in general not
perfectly absorbing.

The second is the use of simultaneous forward and backward sweeps as
opposed to consecutive ones. This idea has been previously tried with
other types of domain decomposition in \cite{Vion2014Thesis}. We find
that this improves the possibilities for parallellization at
very little cost.

%%% some Changes !!!
The third modification is the most interesting from the point of view
of computational cost. We propose to combine the domain decomposition
with a two-grid method, in such a way that the exact inverse at the
coarse level of a two-grid preconditioner is replaced by an
approximate inverse given by a domain decomposition
preconditioner. The result will be called a {\em two-grid
sweeping preconditioner} (TGSP). 
The idea is reminiscent of inner-outer iteration methods. However, we
consider only a single inner iteration. As a result our preconditioner
is a linear map.

Motivating this is the observation that a single iteration of the domain
decomposition preconditioner is considerably more expensive than a
single two-grid or multigrid iteration, compare e.g.\ the computation
times in \cite{PoulsonEtAl2013,CalandraEtAl2013}.  As
a consequence, a single iteration of a TGSP is
considerably cheaper than a single domain decomposition iteration. But
can TGSP also lead to convergence in few iterations? 
Here recent results on multigrid methods for the Helmholtz
equation enter. In \cite{StolkEtAl2014} a class of
multigrid methods for the Helmholtz equation with very good
convergence is studied, based on certain optimized finite difference
discretizations used at the coarse level of the multigrid method.  The
numerical examples below show that the good convergence properties
carry over to the TGSP method, i.e.\ when the exact coarse level
solver is replaced by a domain decomposition preconditioner. The idea
of combining a solver with an outer two grid iteration was previously
studied, using a different setup, in \cite{CalandraEtAl2013}.

A technical complication is the use of multigrid in the presence of
PML layers. This generally requires specifically designed multigrid
methods, e.g.\ in \cite{CalandraEtAl2013} a nonlinear smoother is
used. In this paper we propose two alternatives. The first is the use
of classical absorbing layers, also called sponge layers, instead of
PML layers. The second is a modification in the mesh coarsening in the PML
layers. In this case the use of PML layers of just a few grid cells
wide remains possible. The sponge layers are considerably thicker than
the PML layers, e.g.\ 35 points for the sponge layers in
\cite{RiyantiEtAl2007}, versus around 4 points for the PML layer.

A theoretical result concerning the domain decomposition method with
new transmission and simultaneous forward and backward sweeps is
presented. We show that the method produces an exact solution on the
strip with constant $k$, similar to the domain decomposition method of
\cite{Stolk2013}.

We then study the method using numerical examples. In 2-D we study
problems with up to $7 \cdot 10^{6}$ degrees of freedom, and in 3-D
with up to $10^8$ degrees of freedom. In both cases it is possible to
use quite thin PML layers for the domain decomposition preconditioner,
e.g.\ $w_{\rm pml} = 3$ or $4$ grid cells thick. The convergence of
the method changes very little when simultaneous forward and backward
sweeps are used, compared to executing them after each other.  We show
that for the 3-D examples the two-grid accelerated method indeed leads
to a large reduction in computational cost compared to the ``pure''
sweeping method, and becomes comparable in computation time to the 
fastest methods in the literature.

The setup of the paper is as follows. In section~\ref{sec:dsdd_method}
we describe the double sweep domain decomposition method, including
the modified transmission and simultaneous forward and backward
sweeps.  A theorem describing the behavior of this method on a strip
with constant $k$ is given in section~\ref{sec:dsdd_theorem}.  We then
describe in section~\ref{sec:two_grid} the two-grid sweeping
preconditioner. In section~\ref{sec:implementation} the implementation
will be briefly discussed. Section~\ref{sec:numerical_results}
contains the numerical results. In section~\ref{sec:discussion} a
brief discussion of our results and possible further developments is
given. In an appendix we discuss the discretization of the operators
when PML layers and multigrid are combined using modified mesh
coarsening in the PML layers.

\section{A modified domain decomposition method}
\label{sec:dsdd_method}
In this section we introduce the modified domain decomposition method.

\subsection{Continuous formulation}
\label{subsec:continuous_formulation}
We will formulate the method first in the continuous setting.  We
assume the domain is a rectangle $\Omega = ]0,L[ \times ]0,1[$. It is
straightforward to generalize this to other 2-D and 3-D rectangular domains.

The Helmholtz operator will be denoted by $A$, and is given away from the 
PML or sponge boundary layers by
\begin{equation}
  A = -\partial^2_{xx} - \partial^2_{yy} - k(x,y)^2 .
\end{equation}
In a PML layer at a boundary, say $x=\text{constant}$, it is
obtained by replacing
\begin{equation}
  \frac{\partial}{\partial x} 
\rightarrow
  \frac{1}{1+ i \frac{\sigma_x(x)}{\omega} }
  \frac{\partial}{\partial x} 
\end{equation}
where $\sigma_x = 0$ in the interior of the domain, and positive inside 
the PML layers \cite{Johnson_notespml}. More specifically, motivated
by equation (8) of \cite{Johnson_notespml} we set
\begin{equation}
  \sigma_x = \left\{\begin{array}{ll}
      C_{\rm pml} x^2 & \text{for } x<0\\
      0 & \text{for } 0<x<B_x \\
      C_{\rm pml}(x-B_x)^2 & \text{for } x>B_x
    \end{array}\right .
\end{equation} 
if the PML layers are added outside the domain $x \in [0,B_x]$,
where $C_{\rm pml} = S_{\rm pml} \frac{c_{\rm pml}}{d_{\rm pml}^3}$
with $S_{\rm pml}$ is a dimensionless PML strength parameter, $c_{\rm
  pml}$ is a typical velocity, and $d_{\rm pml}$ the thickness of
the PML layer.
%%% some Changes !!!
In a sponge boundary layer, the constant $k$ is replaced by
$k(1+i\beta(x,y))$. This results in exponential decay of solutions
inside the damping layer, by a factor (in 1-D) of approximately 
$e^{ - \int k \beta(x) \, dx}$. The function $\beta$ was chosen
continuous and quadratically increasing so that in the sponge layer 
a damping on the order of $10^{-2}$ to $10^{-4}$ resulted (note that
reflecting waves pass this layer twice).
Variations in $\beta$ lead to reflections. To make sure that the
reflected energy is small, the sponge layers were several wave lengths
wide. 

Note that absorbing layers of the original domain in general differ
from those introduced in the domain decomposition. In the domain
decomopsition we always use PML layers, of thickness $w_{\rm pml} =
3,4$ or $5$ grid points.  For the original domain we choose between
sponge and PML boundary layers.

We assume the domain is divided in $J$ subdomains 
$] b_{j-1}, b_j [ \times ]0,1[$, with 
\begin{equation}
  0 = b_0 < \ldots < b_J = L ,
\end{equation}
i.e.\ a partition along the $x$-axis.
This partition of the domain will be used for the forward sweep. 
For the backward sweep we assume the domain is divided in $J$ subdomains 
$] \tilde{b}_{j-1}, \tilde{b}_j [ \times ]0,1[$, with 
\begin{equation}
  0 = \tilde{b}_0 < \ldots < \tilde{b}_J = L .
\end{equation}
It is essential that the $b_j$ and the $\tilde{b}_j$ are different and
we will assume that
\begin{equation} \label{eq:ordering_b_btilde1}
  \tilde{b}_j  < b_j < \tilde{b}_{j+1} , \qquad j=1,\ldots,J-1 .
\end{equation}
(A limited number of experiments has been done with $b_j < \tilde{b}_j <
b_{j+1}$, $j=1,\ldots,J-1$ which indicated the method also works well
in this case. Therefore we will formulate the method for both cases.)
Subdomains $\Omega^{(j)}$ (cf.\ equation (10) of \cite{Stolk2013}) are
then defined by
\begin{equation}
  \Omega^{(j)} = ] \min(b_{j-1},\tilde{b}_{j-1}) - L_{\rm PML} (1-\delta_{j,1}),
\max(b_j,\tilde{b_j}) + L_{\rm PML} ( 1- \delta_{j,J}) [        \times ]0,1[ 
\end{equation}
On the domains $\Omega^{(j)}$, functions $k^{(j)}(x,y)$ are defined
that agree with $k$ in the non-PML core of $\Omega^{(j)}$, and are
independent of $x$ and equal to $k$ at the boundary of the core
subdomain inside the added PML layers, i.e.\
\begin{equation}
  k^{(j)}(x,y) = 
\left\{ \begin{array}{ll}
  k(x,y)       & \text{for $\min(b_{j-1},\tilde{b}_{j-1}) \le x \le \max(b_j,\tilde{b_j})$} \\
  k(\min(b_{j-1},\tilde{b}_{j-1}),y) 
              & \text{for $x < \min(b_{j-1},\tilde{b}_{j-1})$ (if $j>1$)} \\
  k(\max(b_j,\tilde{b_j}),y)    
              & \text{for $x > \max(b_j,\tilde{b_j})$ (if $j<J$).}
\end{array} \right.
\end{equation}
On the domains $\Omega^{(j)}$ operators $A^{(j)}$ are defined as Helmholtz
operators with PML modifications, similar as $A$ was defined on $\Omega$.

To derive the method for transmission, we consider the case
$J=2$. Then, in the forward sweep, the equation is first solved on
$\Omega^{(1)}$ with $f^{(1)} = H(b_1-x) f$ as right hand side, where
$H$ denote the Heaviside function.  Subsequently it is solved on
$\Omega^{(2)}$ with as right hand side $f^{(2)} = H(x-b_1) f$ plus a
contribution from the local solution on $\Omega^{(1)}$, which is to be
determined. Suppose $v^{(1)}$ is the solution of $A^{(1)} v^{(1)} =
f^{(1)}$. Then, ideally we would like to obtain $w$ such that
$H(b_1-x) v^{(1)} + w$ is the true solution, in other words
\begin{equation} \label{eq:idea_new_transmission_eq1}
  A ( H(b_1-x) v^{(1)} + w ) = f 
\end{equation}
(cf.\ \cite{SchadleEtAl2007}).
Then $w$ must satisfy
\begin{equation} \label{eq:idea_new_transmission}
  A w = f - A (  H(b_1-x) v^{(1)} )
  = f^{(2)} - A^{(1)}  (  H(b_1-x) v^{(1)} )  +  H(b_1-x) A^{(1)} v^{(1)} .
\end{equation}
To arrive at a domain decomposition method, we observe that the right
hand side is supported in the set $x \in [b_1,b_2]$ and solve this on
$\Omega^{(2)}$, i.e.\ we solve
\begin{equation} \label{eq:idea_new_transmission_eq3}
  A^{(2)} v^{(2)} 
  = f^{(2)} - A^{(1)}  (  H(b_1-x) v^{(1)} )  +  H(b_1-x) A^{(1)} v^{(1)} 
\end{equation}
The second and third terms on the right hand side amount to the
transmission of information from the solution on subdomain 1 to the
equation for subdomain 2. Below we will show that they generate a
forward propagating wave in the subdomain $x \in [b_1,b_2]$, thereby
extending the truncated solution $H(b_1-x)v^{(1)}$.  We set $u =
v^{(1)} + v^{(2)}$ as approximate solution. We show below this can
model forward propagating waves over the entire domain, but not the
backward propagating waves. These can be computed in a backward sweep:
solving first on subdomain 2 and then on subdomain 1. Waves reflecting
back and forth between the subdomains can be obtained in an
iterative process.

To denote the contribution from neighboring solutions to the right
hand sides of some subdomain, we will define transmission operators
$T^{(j)}$ and $\tilde{T}^{(j)}$, for the forward and backward sweep
respectively. The operator $T^{(j)}$ acts on $v^{(j-1)}$ defined on
$\Omega^{(j-1)}$ (where $v^{(j-1)}$ must be such that the product
$H(b_{j-1} - x) v^{(j-1)}$ is well defined), and is defined by
\begin{equation} \label{eq:define_transmission_Ttilde}
  T^{(j)} v^{(j-1)} = - A^{(j-1)} ( H(b_{j-1} - x) v^{(j-1)} 
    + H(b_{j-1} - x) A^{(j-1)} v^{(j-1)} .
\end{equation}
This is a distribution supported on $x = b_{j-1}$ and hence can be
considered a distribution on $\Omega^{(j)}$. 
Similarly we define $\tilde{T}^{(j)}$ by
\begin{equation}
  \tilde{T}^{(j)} w^{(j+1)} = - A^{(j+1)} ( H(x - \tilde{b}_j) w^{(j+1)} 
    + H(x-\tilde{b}_j) A^{(j+1)} w^{(j+1)} .
\end{equation}

We can now describe the domain decomposition method with the forward
and backward sweeps performed after each other. By $I_{x \in [\alpha,\beta]}$
we denote the indicator function which is one for $x \in
[\alpha,\beta]$ and we will assume $T^{(1)} = \tilde{T}^{(J)} = 0$.
The domain decomposition preconditioner is then described by 
the algorithm \proc{SweepingPrecUDContinuous} in Table~\ref{tab:algorithms1and2}.
\begin{table}
%the following algorithm:
%\begin{equation} \label{eq:UDsweep_algorithm_continuous}
%\begin{minipage}[c]{12cm}
\begin{codebox}
\Procname{\proc{SweepingPrecUDContinuous}$(f)$}
\li $u=0$
\li \For $j = 1, \ldots, J$
\li     \Do
            solve $v^{(j)}$ from
            $A^{(j)} v^{(j)} = I_{x \in [b_{j-1},b_j]} f + T^{(j)} v^{(j-1)}$
\li         $u = u + I_{x \in [b_{j-1},b_j]} v^{(j)}$
        \End
\li $g = f - A u$ 
\li \For $j=J,J-1,\ldots,1$
\li     \Do
            solve $w^{(j)}$ from 
            $A^{(j)} w^{(j)} = I_{x \in [\tilde{b}_{j-1},\tilde{b}_j]}
            g + \tilde{T}^{(j)} w^{(j+1)}$
\li         $u = u + I_{x \in [\tilde{b}_{j-1},\tilde{b}_j]} w^{(j)}$
        \End
\li \Return u
\end{codebox}
\begin{codebox}
\Procname{\proc{SweepingPrecXContinuous}$(f)$}
\li $u=0$
\li \For $j = 1, \ldots, J/2$
\li     \Do
            solve $v^{(j)}$ from
            $A^{(j)} v^{(j)} = I_{x \in [b_{j-1},b_j]} f + T^{(j)} v^{(j-1)}$
\li         $u = u + I_{x \in [b_{j-1},b_j]} v^{(j)}$
\li         \If $j \neq 1$
\li             \Then
                    $k=J+2-j$
\li                 solve $v^{(k)}$ from
                        $A^{(k)} v^{(k)} = I_{x \in [\tilde{b}_{k-1},\tilde{b}_k]} f 
                        + \tilde{T}^{(k)} v^{(k+1)}$
\li                 $u = u + I_{x \in [\tilde{b}_{k-1},\tilde{b}_k]} v^{(k)}$
                \End
        \End
\li $j=J/2+1$
\li solve $v^{(j)}$ from $A^{(j)} v^{(j)} = I_{x \in [b_{j-1},\tilde{b}_j]} f 
        + T^{(j)} v^{(j-1)} + \tilde{T}^{(j)} v^{(j+1)}$
\li $u = u + I_{x \in [b_{j-1},\tilde{b}_j]} v^{(j)}$
\li $g = f - A u$ 
\li solve $w^{(j)}$ from $A^{(j)} w^{(j)} = I_{x \in [\tilde{b}_{j-1},b_j]} g$
\li $u = u + I_{x \in [\tilde{b}_{j-1},b_j]} w^{(j)}$
\li \For $j=J/2, J/2-1, \ldots, 1$
\li     \Do
            solve $w^{(j)}$ from 
            $A^{(j)} w^{(j)} = I_{x \in [\tilde{b}_{j-1},\tilde{b}_j]} g + \tilde{T}^{(j)} w^{(j+1)}$
\li         $u = u + I_{x \in [\tilde{b}_{j-1},\tilde{b}_j]} w^{(j)}$
\li         \If $j < J/2$
\li             \Then
                    $k=N+1-j$
\li                 solve $w^{(k)}$ from
                        $A^{(k)} w^{(k)} = I_{x \in [b_{k-1},b_k]} g  
                        + T^{(k)} w^{(k-1)}$
\li                 $u = u + I_{x \in [b_{k-1},b_k]} w^{(k)}$
                \End
        \End
\li \Return u
\end{codebox}
\caption{Domain decomposition algorithms in the continuous setting}
\label{tab:algorithms1and2}
\end{table}
%\end{minipage}
%\end{equation}

Note that the restrictions of $g$ to the subdomains $x \in
[\tilde{b}_{j-1},\tilde{b}_j]$ are well defined, because the singular
support of $v$ is at the boundaries $x = b_j$. Similarly, the singular
support of the residual $f - A u$ is at the boundaries $x =
\tilde{b}_j$, so that, in the next iteration of a preconditioned
iterative solver, the restrictions of the residual $f - Au$ to the
sets $x \in [b_{j-1},b_j]$ are well defined.

Next we consider the continuous formulation of a domain decomposition
method with {\em simultaneous sweeps}. We will also refer to this as
intersecting sweeps or X-sweep, because, in a plot of the subdomain 
being solved versus the step number in the algorithm, the resulting 
graph contains two intersecting lines like a diagonal cross.
We assume that $J$ is even and that this intersection is at a particular subdomain numbered
$j_{\rm mid}$, chosen such that $j_{\rm mid} = J/2  +1$.  
The algorithm for this case is algorithm \proc{SweepingPrecXContinuous} in
Table~\ref{tab:algorithms1and2}.
Just like above, the restrictions of $g$ to the subdomains are well
defined because the $b_j$ are different from the $\tilde{b}_j$.

The resulting solutions $u$ for the algorithms in
Table~\ref{tab:algorithms1and2} depend linearly on $f$ and will be
denoted by $P_{\rm UD} f$ and $P_{\rm X} f$ respectively.

\subsection{Discrete formulation}

For the discrete formulation we assume that $A$ is discretized on a
regular or rectilinear mesh. The mesh is to consist of $N_x \times
N_y$ cells. Because we use Dirichlet boundary conditions, there are
$(N_x - 1) \times (N_y-1)$ unknowns. If we denote the degrees of
freedom by $u_{i,j}$, we will write the discretized Helmholtz equation
as 
\begin{equation}
  (A u)_{i,k} = \sum_{\tilde{i},\tilde{k}} a_{i,k;\tilde{i},\tilde{k}}
  u_{\tilde{i},\tilde{k}} ,
\end{equation}
We will assume a compact stencil discretization, i.e.\
$a_{i,k;\tilde{i},\tilde{k}} = 0$ if $|i-\tilde{i}| > 1$ or $|k -
\tilde{k}| > 1$.

The subdomain boundaries $b_j$ and $\tilde{b}_j$ are assumed to be at
half grid points $x_{\beta_j + 1/2}$ and
$x_{\tilde{\beta}_j+1/2}$. The discrete equivalent to the interval
$]b_{j-1},b_j[$ is therefore the set of points $\{ x_{\beta_{j-1}+1} ,
\ldots, x_{\beta_j} \}$. After the first set of discrete subdomain 
boundaries $\beta_j$ is chosen, the second set is defined by
\begin{equation}
\begin{split}
  \tilde{\beta}_0 = {}& \beta_0
\\
  \tilde{\beta}_J = {}& \beta_J
\\
  \tilde{\beta}_j = {}& \beta_j - 1
\qquad  \text{for  $j=1,\ldots, J-1$}.
\end{split}
\end{equation}

The discretized transmission matrix $T^{(j)}$ is a matrix from the layers
with global coordinates $i = \beta_{j-1}, \beta_{j-1}+1$ in subdomain
$j-1$, to the layers with the same global coordinates in subdomain
$j$.  We define operators $J_{\rm out}^{(j-1)}$ to extract these layers
from the unknown on subdomain $j-1$, and operators $( J_{\rm in}^{(j)}
)^t$ to inject (is adjoint of restriction) into subdomain $j$. It is
straightforward to show that the discretized transmission operator,
defined using (\ref{eq:define_transmission_Ttilde}), is then given by a product
\begin{equation}
  (J_{\rm in}^{(j)})^t  T^{(j)}  J_{\rm out}^{(j-1)}
\end{equation}
where the discrete operator $T^{(j)}$ is given by
(note that $s,\tilde{s} \in \{ 0,1 \}$)
\begin{equation}
  T^{(j)}_{1+s,k;1+\tilde{s},\tilde{k}}
  = \left\{\begin{array}{ll}
    0 & \text{when $\tilde{s} = s$} \\
    \pm A_{\beta_{j-1}+s,k;\beta_{j-1}+\tilde{s},\tilde{k}} &
        \text{when $\tilde{s} - s = \pm 1$,}
        \end{array} \right. .
\end{equation}
(We use the same notation $T^{(j)}$ for the continuous and discrete
transmission operators, from the context it should be clear which one is
intended.)
Let operators $\tilde{J}_{\rm out}^{(j+1)}$ to extract these layers
from the unknown on subdomain $j+1$, and operators $(\tilde{J}_{\rm
  in}^{(j)})^t$ be defined similarly to extract layers with global
coordinates $i = \tilde{\beta}_j,\tilde{\beta}_j+1$ from subdomain
$j+1$, and to inject them into subdomain $j$.
The discrete transmission matrix in this case has components
\begin{equation}
  \tilde{T}^{(j)}_{1+s,j;1+\tilde{s},\tilde{j}}
  = \left\{\begin{array}{ll}
    0 & \text{when $\tilde{s} = s$} \\
    \mp A_{\tilde{\beta}_j+s,j;\tilde{\beta}_j+\tilde{s},\tilde{j}} &
        \text{when $\tilde{s} - s = \pm 1$}
        \end{array} \right.
\end{equation}
(again $s,\tilde{s} \in \{ 0,1 \}$).

To map data between subdomains and the full domain we define
$J(j,a,b)$ to be the matrix that maps degrees of freedom $u_{i,k}$ with
$i \in \{ a+1, \ldots, b \}$ to the corresponding degrees of freedom
for a discrete function defined on $\Omega^{(j)}$. The transpose
$J(j,a,b)^T$ maps values from the a discrete function on the subdomain
to a discrete function of the full domain.

With these definitions and results we can define algorithms for the
discrete domain decomposition preconditioners that were presented
above in the continuous setting. A few helper algorithms are presented
in Table~\ref{tab:sweeping_prec_helper_algorithms}.  The algorithm
$\proc{SubdomSolve}$ performs a generic subdomain solve and update
including the handling of transmission data. The argument $j$ is the
subdomain number; $a,b$ describe which layers of degrees of freedom
are to be copied from the right hand side on $\Omega$ to the right
hand side on $\Omega^{(j)}$; $\tilde{a}$, $\tilde{b}$ describe which
layers from to solution on $\Omega^{(j)}$ to copy to the approximate
solution on $\Omega$; flags $\tau_j$, $j=1,2,3,4$ indicate whether
transmission is done for (in,forward), (out,forward), (in,backward)
and (out,backward) uses respectively and the $B_j$ are variables used
for storing or retrieving transmission data.  The algorithms
$\proc{ForwardSweep}$ and $\proc{BackwardSweep}$ execute a series of
solves, using the transmission matrices. They have as arguments the
right hand side and unknown for the approximate solution, the first
and last subdomain to be included and a buffer $B$ to store
transmission data.

The preconditioner applications, using non-simultaneous and
simultaneous forward and backward sweeps are given in
Table~\ref{tab:sweeping_prec_main_algorithms}.  We have included an
algorithm for domain decomposition with partial sweeps called
$\proc{SweepingPrecNX}$. In \cite{Vion2014Thesis} such an algorithm
was given for domain decomposition with different
interface/transmission conditions, the equivalent for our method is
including in Table~\ref{tab:sweeping_prec_main_algorithms}.  In this
algorithm intersecting sweeps are done over groups of subdomains. The
boundary domains of these groups are given by $j_{\rm cell,m}$, 
$m=0, \ldots, N_{\rm cell}$, and the center domains where the local 
sweeps intersect are given by $j_{\rm mid,m}$, 
$m=1, \ldots, N_{\rm cell}$. It is assumed that $j_{\rm cell,0} = -1$ and 
$j_{\rm cell,N_{\rm cell}} = J + 1$.

\begin{table}
\begin{codebox}
\Procname{$\proc{SubdomSolve}(u,f,j,a,b,\tilde{a},\tilde{b},
    \tau_1,B_1,\tau_2,B_2,\tau_3,B_3,\tau_4,B_4)$}
\li $f_{\rm d} = J(j,a,b) f$
\li \If $\tau_1$
        \Then
\li         $f_{\rm d} = f_{\rm d} + (J_{in}^{(j)})^y T^{(j)} B_1$
        \End
\li \If $\tau_3$
        \Then
\li         $f_{\rm d} = f_{\rm d} + (\tilde{J}_{\rm in}^{(j)})^t \tilde{T}^{(j)} B_3$
        \End
\li $u_{\rm d} = (A^{(j)})^{-1} f_{\rm d}$
\li $u = u + J(j,\tilde{a},\tilde{b})^T u_{\rm d}$
\li \If $\tau_2$
        \Then
\li         $B_2 = J_{\rm out}^{(j)} u_{\rm d}$
        \End
\li \If $\tau_4$
        \Then
\li         $B_4 = \tilde{J}_{\rm out}^{(j)} u_{\rm d}$
        \End
\end{codebox}
\begin{codebox}
\Procname{$\proc{ForwardSweep}(u,f,j_0,j_1,B)$}
\li \For $j = j_0, j_0+1, \ldots, j_1$
\li     \Do
            $\proc{SubdomSolve}(u,f,j,\beta_{j-1},\beta_j,\beta_{j-1},\beta_j,
                 j>1,B,j<N_{\rm dom},B,0,\NIL,0,\NIL)$
        \End
\end{codebox}
\begin{codebox}
\Procname{$\proc{BackwardSweep}(u,f,j_0,j_1,B)$}
\li \For $j = j_0,j_0-1, \ldots, j_1$
\li     \Do
            $\proc{SubdomSolve}(u,f,j,
                   \tilde{\beta}_{j-1},\tilde{\beta}_j,\tilde{\beta}_{j-1},\tilde{\beta}_j,
                   0,\NIL,0,\NIL,j<N_{\rm dom},B,j>1,B)$
        \End
\end{codebox}
\begin{codebox}
\Procname{$\proc{MidSolveIn}(u,f,j,B_1,B_2)$}
\li $\proc{SubdomSolve}(u,f,j,\beta_{j-1},\tilde{\beta_j},\beta_{j-1},\tilde{\beta_j},
    1,B_1,0,\NIL,1,B_2,0,\NIL)$
\end{codebox}
\begin{codebox}
\Procname{$\proc{MidSolveOut}(u,f,j,B_1,B_2)$}
\li $\proc{SubdomSolve}(u,f,j,\tilde{\beta}_{j-1},\beta_j,\tilde{\beta}_{j-1},\beta_j,
    0,\NIL,1,B_1,0,\NIL,1,B_2)$
\end{codebox}
\caption{Helper algorithms}
\label{tab:sweeping_prec_helper_algorithms}
\end{table}

\begin{table}
\begin{codebox}
\Procname{$\proc{SweepingPrecUD}(f)$}
\li $u \gets 0$
\li $\proc{ForwardSweep}(u,f,1,N_{\rm dom},B)$
\li $g = f - A u$
\li $\proc{BackwardSweep}(u,g,N_{\rm dom},1,B)$
\li \Return $u$ 
\end{codebox}
\begin{codebox}
\Procname{$\proc{SweepingPrecX}(f)$}
\li $u \gets 0$
\li $\proc{ForwardSweep}(u,f,1,j_{\rm mid}-1,B_1)$
\li $\proc{BackwardSweep}(u,f,J,j_{\rm mid}+1,B_2)$
\li $\proc{MidSolveIn}(u,f,j_{\rm mid},B_1,B_2)$
\li $g = f - A u$
\li $\proc{MidSolveOut}(u,g,j_{\rm mid},B_2,B_1)$
\li $\proc{BackwardSweep}(u,g,j_{\rm mid}-1,1,B_1)$
\li $\proc{ForwardSweep}(u,g,j_{\rm mid}+1,J,B_2)$
\li \Return $u$ 
\end{codebox}
\begin{codebox}
\Procname{$\proc{SweepingPrecNX}(f)$}
\li \For $m=1,\ldots,N_{\rm cell}$
        \Do
\li         \If $m<N_{\rm cell}$
                \Then
\li                 $\proc{MidSolveOut}(u,f,j_{{\rm cell},m},B_{2m-1})$
                \End
\li         $\proc{ForwardSweep}(u,f,j_{{\rm cell},m-1}+1,j_{{\rm mid},m}-1,B_{2m-1})$
\li         $\proc{BackwardSweep}(u,f,j_{{\rm cell},m}-1,j_{{\rm mid},m}+1,B_{2m})$
\li         $\proc{MidSolveIn}(u,f,j_{{\rm mid},m},B_{2m-1},B_{2m})$
        \End
\li $g = f - A u$
\li \For $m=1,\ldots,N_{\rm cell}$
        \Do
\li         $\proc{MidSolveOut}(u,g,j_{{\rm mid},m},B_{2m},B_{2m-1})$
\li         $\proc{BackwardSweep}      (u,g,j_{{\rm mid},m}-1,j_{{\rm cell},m-1}+1,B_{2m-1})$
\li         $\proc{ForwardSweep}        (u,g,j_{{\rm mid},m}+1,j_{{\rm cell},m}-1,B_{2m})$
\li         \If $m<N_{\rm cell}$
                \Then
\li                 $\proc{MidSolveIn}(u,g,j_{{\rm cell},m},B_{2m},B_{2m+1})$
                \End
        \End
\li \Return u
\end{codebox}
\caption{Algorithms $\proc{SweepingPrecUD}$, $\proc{SweepingPrecX}$ and
$\proc{SweepingPrecNX}$ for different variants of the sweeping preconditioner.}
\label{tab:sweeping_prec_main_algorithms}
\end{table}

\section{Theoretical results}
\label{sec:dsdd_theorem}

Here we study the domain decomposition in case of constant $k$ on a
line segment in one dimension and for a two-dimensional strip with PML
layers only at the boundaries $x=0$ and $x=L$.

\subsection{One-dimensional analysis}
\label{subsec:one_dimensional_analysis}
We will show that the domain decomposition method reproduces the exact
solution when the domain is a line segment and $k$ is constant.

In one dimension absorbing boundary conditions are given by
Robin boundary conditions and the problem on $]0,L[$ becomes
\begin{equation} \label{eq:1D_problem}
\begin{split}
-\partial_{xx}^2 u(x) - k^2 u(x) = {}& f(x)
\\
  \partial_x u(0) + iku(0) = {}& h_1 ,
\\
- \partial_x u(L) + iku(L) = {}& h_2  
\end{split}
\end{equation}
One can also enlarge the domain, i.e.\ if $\alpha \le 0 < L \le \beta$
the problem on can be considered on $]\alpha,\beta[$ with boundary
conditions at $\alpha$,$\beta$, without affecting the solution on
$]0,L[$, because in each case an unbounded domain is simulated.
The solution for (\ref{eq:1D_problem}) is given by
\begin{equation} \label{eq:solutionformula1D}
  u(x) = \frac{i}{2k} \int_{0}^x e^{ik(x-s)} f(s) \, ds 
+ \frac{i}{2k} \int_x^L e^{-ik(x-s)} f(s) \, ds
+ \frac{e^{ikx}}{2ik} h_1
+ \frac{e^{-ik(x-L)}}{2ik} h_2 
\end{equation}

In some case we are interested in solutions $w$ to the
\begin{equation} \label{eq:diffeq_with_minusAu}
  A w = f - A u
\end{equation}
on an interval $]\alpha,\beta[$ with homogeneous boundary conditions
$\partial_x w(\alpha) + ikw(\alpha) = 0$ 
and $- \partial_x w(\beta) + ik w(\beta) = 0$.
In this case we determine
\begin{equation}
\begin{split}
  R_+ u(\alpha) \stackrel{\rm def}{=} {}& 
    \frac{1}{2ik} \left( \partial_x u(\alpha) + iku(\alpha) \right)
\\
  R_- u(\beta) \stackrel{\rm def}{=} {}& 
    \frac{1}{2ik} \left( -\partial_x u(\beta) + iku(\beta) \right)
\end{split}
\end{equation}  
The solution to (\ref{eq:diffeq_with_minusAu}) then satisfies
\begin{equation} \label{eq:sol_with_minusAu}
\begin{split}
  u(x) + w(x) = {}&
\frac{i}{2k} \int_{\alpha}^x e^{ik(x-s)} f(s) \, ds 
+ \frac{i}{2k} \int_x^\beta e^{-ik(x-s)} f(s) \, ds 
\\
  {}& +e^{ik(x-\alpha)} R_+u(\alpha) + e^{-ik(x-\beta)} R_-u(\beta)
\end{split}
\end{equation}

The effect of using a transmission source $T^{(j)} v^{(j-1)}$ can be
analyzed using equations (\ref{eq:diffeq_with_minusAu}) to
(\ref{eq:sol_with_minusAu}). We will consider the case $J=2$ given in
equations (\ref{eq:idea_new_transmission_eq1}) to
(\ref{eq:idea_new_transmission_eq3}). First note that for $0 < x <
b_1$
\begin{equation}
  v^{(1)} = \frac{i}{2k} \int_0^x e^{ik(x-s)} f(s) \, ds
    + \frac{i}{2k} \int_x^{b_1} e^{-ik(x-s)} f(s) \, ds .
\end{equation}
Using that we can enlarge the domain, we consider equation 
(\ref{eq:idea_new_transmission_eq3}) as an equation of the
type (\ref{eq:diffeq_with_minusAu}) on the interval
$x \in ]b_1-\epsilon,b_2[$. This gives that for 
$x \in ]b_1-\epsilon,b_2[$ we have
\begin{equation}
\begin{split}
  v^{(2)}(x) + v^{(1)} H(b_1 - x)
  = {}& e^{ik(x - (b_1-\epsilon))} R_+ v^{(1)}(b_1 - \epsilon)
\\
{}& + \frac{i}{2k} \int_{b_1 - \epsilon}^x e^{ik(x-s)} f(s) \, ds 
    + \frac{i}{2k} \int_x^L e^{-ik(x-s)} f(s) \, ds .
\end{split}
\end{equation}
Considering that $v^{(1)} H(b_1 - x) = 0$ for $x>b_1$, that 
$R_+ v^{(1)}(b_1) = \frac{i}{2k} \int_0^{b_1} e^{ik(b_1-s)} f(s) \, ds$,
and taking the limit $\epsilon \rightarrow 0$, we obtain for 
$x \in ]b_1,b_2[$
\begin{equation}
  v^{(2)}(x) = \frac{i}{2k} \int_0^x e^{ik(x-s)} f(s) \, ds 
    + \frac{i}{2k} \int_x^{b_2} e^{-ik(x-s)} f(s) ,
\end{equation}
which is the correct solution on this subdomain. Similarly it follows
that the effect of the transmission source $T^{(j)} v^{(j-1)}$ in the
right hand side of 
$A^{(j)} v^{(j)} = I_{x \in [b_{j-1},b_j]} f + T^{(j)} v^{(j-1)}$ is a 
contribution
\begin{equation}
  e^{ik(x-b_{j-1})} R_+v^{(j-1)}(b_{j-1}) 
\end{equation}
to the solution $v^{(j)}$ on $]b_{j-1},b_j[$.

By induction we then find the following for the forward sweep in
algorithm \proc{SweepingPrecUDContinuous}. After step $j$ in the loop,
we have 
\begin{equation} \label{eq:induction_forward_sweep}
  u(x) = \frac{i}{2k} \int_0^x e^{ik(x-s)} f(s) \, ds 
+ \frac{i}{2k} \int_x^{b_l} e^{-ik(x-s)} f(s) \, ds ,
\qquad 
  \text{for $x \in ]b_{l-1},b_l[$, $l\le j$.}
\end{equation}
and $u(x) = 0$ for $x > b_j$.
For the backward sweep (\ref{eq:sol_with_minusAu}) is used again.
By induction one can show that after subdomain $j$ is updated, the
solution is given by
\begin{equation}
  u(x) = \frac{i}{2k} \int_0^x e^{ik(x-s)} f(s) \, ds 
+ \frac{i}{2k} \int_x^L e^{-ik(x-s)} f(s) \, ds ,
\end{equation}
for $x > b_{j-1}$ while $u(x)$ is still given by 
(\ref{eq:induction_forward_sweep}) for $x \in [b_{l-1},b_l]$, $l<j$.
Hence algorithm \proc{SweepingPrecUDContinuous} yields the correct
solution. 

For the simultaneous sweeps, similarly after step $j$ of the first
loop we have
\begin{equation}
  u(x) = \left\{ \begin{array}{ll}
\frac{i}{2k} \int_0^x e^{ik(x-s)} f(s) \, ds 
+ \frac{i}{2k} \int_x^{b_l} e^{-ik(x-s)} f(s) \, ds 
& \text{for $x \in [b_{l-1},b_l]$, $l\le j$}
\\
\frac{i}{2k} \int_{\tilde{b}_{l-1}}^x e^{ik(x-s)} f(s) \, ds 
+ \frac{i}{2k} \int_x^L e^{-ik(x-s)} f(s) \, ds 
& \text{for $x \in [\tilde{b}_{l-1},\tilde{b}_l]$, $l \ge J+2-j$} 
\\
0 & \text{otherwise} . 
\end{array} \right.
\end{equation}
After lines 9-11 of the algorithm the function $u$ satisfies for 
$x \in ]b_{J/2},\tilde{b}_{J/2+1}[$ 
\begin{equation}
  u(x) = \int_0^x e^{ik(x-s)} f(s) \, ds 
+ \frac{i}{2k} \int_x^L e^{-ik(x-s)} f(s) \, ds ,
\end{equation}
which is the true solution. Next one can show inductively that 
steps 13-21 in the algorithm yield the correct solution in each
subdomain that is updated, implying that the algorithm 
\proc{SweepingPrecXContinuous} yields the correct solution.

\subsection{Modified domain decomposition method on the strip}

We next consider the problem with $k = \text{constant}$ on the strip
$]0,L[ \times ]0,1[$, with Dirichlet boundary conditions at $y=0$ and
$y=1$ and PML boundary layers at $x=0$ and $x=L$. In this section we
will assume that a PML layer behaves like a perfect non-reflecting
boundary condition. In essence we will show that Theorem 1 of
\cite{Stolk2013} remains valid for the modified method.

After a Fourier transform in $y$ the solution becomes of the form
$u = \sum_l \sin(2 \pi l) \hat{u}_l(x)$,
$l=1,2,\ldots$, and writing $\hat{u}_l(x) =\hat{u}(x,\eta)$, $\eta =
2\pi l$, the Helmholtz equation becomes a family of ODE's that reads
\begin{equation} \label{eq:Fourier_transformed_PDE_2D}
  -\partial_{xx}^2 \hat{u} + \eta^2 \hat{u} - k^2 \hat{u} 
    = \hat{f}(x,\eta)
\end{equation}
We assume that $k \neq 2\pi l$ for all integers $l>0$.
The non-reflecting boundary condition becomes
\begin{align} \label{eq:perfectly_non_reflecting_bc1}
   \partial_x \hat{u} + \lambda \hat{u} = {}& h_1
&&
  \text{at $x=0$}
\\ \label{eq:perfectly_non_reflecting_bc2}
 - \partial_x \hat{u} + \lambda \hat{u} = {}& h_2
&&
  \text{at $x=L$,}
\end{align}
where $\lambda$ is given by
\begin{equation} \label{eq:perfectly_non_reflecting2}
 \lambda 
= \left\{ \begin{array}{ll}
  i \sqrt{k^2 - \eta^2} & \text{if } |\eta| < k \\
  - \sqrt{\eta^2 - k^2} & \text{if } |\eta| > k ,
\end{array}\right.
\end{equation}
and $h_1$ and $h_2$ are 0 for homogeneous non-reflecting boundary
conditions and non-zero if incoming waves are to be modeled.

In this case we can apply exactly the same analysis as in
section~\ref{subsec:one_dimensional_analysis} to the problems for each
$\eta$. For example, the solution formula for
(\ref{eq:Fourier_transformed_PDE_2D}-\ref{eq:perfectly_non_reflecting_bc2})
is straightforwardly derived and given by
\begin{equation}
\hat{u}(x,\eta) =
\frac{-1}{2\lambda} \int_0^x e^{\lambda (x-s)} \hat{f}(s,\eta) \, ds
  + \frac{-1}{2\lambda} \int_x^L e^{-\lambda(x-s)} \hat{f}(s,\eta) \, ds
  + \frac{e^{\lambda x}}{2\lambda} h_1
  + \frac{e^{-\lambda (x-L)}}{2\lambda} h_2
\end{equation}
Thus we have

\begin{theorem}
On the strip $]0,L[ \times ]0,1[$ with absorbing boundaries at $x=0$
and $x=L$ and constant $k$, the map $P_{\rm X}$ satisfies $A P_{\rm X} f = f$.
\end{theorem}

\section{Two-grid domain decomposition preconditioner}
\label{sec:two_grid}

In this section we describe a method in which a domain decomposition
preconditioner is used as an inexact coarse level solver in a two-grid
method. We consider the case where a two-grid cycle is used as
preconditioner for GMRES. The modified two-grid cycle, with domain
decomposition preconditioner used as coarse level solver, will be
called a two-grid sweeping preconditioner or TGSP. It follows from
computation times given in \cite{PoulsonEtAl2013,CalandraEtAl2013}
that a TGSP application is considerably cheaper than a direct sweeping
preconditioner application. Since the cost of a solve is roughly given
by the cost of a preconditioner application times the number of
iterations, the question is what happens with the number of iterations
when an outer two-grid iteration is added.

In \cite{StolkEtAl2014} it is shown that a certain class of two-grid
methods converges rapidly. This of course refers to the case using an
exact coarse level solver. A priori it is unknown whether these good
convergence properties extend to the case of an inexact, domain
decomposition based coarse level solver, also because in the multigrid
method the sweeping preconditioner is applied at coarser meshes than
it has been tested with so far, using e.g.\ five instead of ten points
per wavelength.  However, it is clear that an efficient solver would
result if the convergence doesn't degrade too much. 

The purpose of the present section is to describe a two-grid sweeping
preconditioner based on the two-grid method of
\cite{StolkEtAl2014}. In sections below we will show that in numerical
examples the convergence remains good and that the method is in fact
highly efficient.

In two subsections we will separately discuss the cases with and
without PML boundary layers present. The presence of PML layers makes
it necessary to modify the multigrid method. We opt for a specific
modification where the mesh coarsening in the PML layers is
changed. Alternatively the smoother can be modified, see e.g.\
\cite{CalandraEtAl2013}. When PML layers are absent we use classical
damping layers as absorbing layers near the boundary of the domain
$\Omega$. See \cite{TrottenbergOosterleeSchueller2001} for background
on multigrid methods.

The original problem will be standard second order finite
differences. The discretization on a regular mesh of the 1-D second
order operator $u \mapsto - \pdpd{}{x} \left( \beta(x) \pdpd{u}{x}
\right)$ is given by
\begin{equation}
  h^{-2} \left( - \beta_{i-1/2} u_{i-1} 
  + (\beta_{i-1/2} + \beta_{i+1/2}) u_i
  - \beta_{i+1/2} u_{i+1} \right) .
\end{equation}
This formula is used to find the following 5-pt finite difference
discretization of the Helmholtz equation (in 2-D) in presence of PML
boundary layers
\begin{equation} \label{eq:define_pml_fd}
\begin{split}
& \frac{1}{h^2 \alpha_{2,k}}
  \left( - \alpha_{1,i-1/2} u_{i-1,k} 
  + (\alpha_{1,i-1/2} + \alpha_{1,i+1/2}) u_{i,k}
  - \alpha_{1,i+1/2} u_{i+1,k} \right) 
\\
& + \frac{1}{h^2 \alpha_{1,i}}
  \left( - \alpha_{2,k-1/2} u_{i,k-1} 
  + (\alpha_{2,k-1/2} + \alpha_{2,k+1/2}) u_{i,k}
  - \alpha_{2,k+1/2} u_{i,k+1} \right) 
\\
& - \frac{k_{i,k}^2}{\alpha_{1,i} \alpha_{2,k}}
  = \frac{1}{\alpha_{1,i} \alpha_{2,k}} f_{i,k} ,
\end{split}
\end{equation}
where where $\alpha_j(x_j) = \frac{1}{1 + i \omega^{-1}
  \sigma_j(x_j)}$ (with $j=1,2$ referring to the $x$ and $y$ axes
respectively). (In absence of PML boundary layers, the coefficients
$\alpha_j$ are equal to $1$.)

\subsection{The two-grid method in absence of PML layers}

In this subsection we will discuss the two-grid method to be used in
absence of PML boundary layers.  This method is according to
\cite{StolkEtAl2014}. It is based on the V-cycle, full weighting
prolongation and restriction operators and $\omega$-Jacobi smoothers,
with parameters given in section~\ref{sec:numerical_results} below.
As mentioned, the two-grid method is used as preconditioner for GMRES.

The main difference of the method of \cite{StolkEtAl2014} compared to
standard multigrid methods is that optimized finite difference
operators constructed in that paper are used as coarse level
discretization.  These are designed such that phase speed differences
between fine and coarse level discretizations are minimal. We recall
the definition of these operators in appendix~\ref{sec:optimized_fd}
that treats coarse level discretizations for the case that PML layers
are present. A second difference is in the choice of parameters for
the smoother. In order to have good convergence the weight $\omega$ in
the $\omega$-Jacobi smoother and the number of pre- and postsmoothing
steps $\nu$ used in the V-cycle must be chosen quite
specifically. Results in \cite{StolkEtAl2014} show that convergence
properties depend sensitively on these parameters.

The inclusion of an inexact, domain decomposition based coarse level
solver is done straightforwardly: The coarse level solver is simply
replaced by a preconditioner application. This is of course an additional
difference with standard multigrid. The parameters (number of
subdomains, PML width and PML strength $S_{\rm pml}$) will be
discussed below in the section on numerical examples.

\subsection{Using PML layers in the two-grid method}
\label{sec:two_grid_pml}

With PML-layers it is typically necessary to modify the multigrid
method because convergence becomes poor. It is not easy to precisely
pinpoint the cause of this behavior. The local Fourier analysis of the
Helmholtz operator without PML is inapplicable for two reasons. First
the matrix is changed locally, and second the coefficients $\sigma_x$,
and $\sigma_y$ vary rapidly, implying that the assumptions of the
local Fourier analysis are not valid. These are also the potential
reasons for which convergence is hampered.

\begin{figure}
\begin{center}
\includegraphics[height=28mm]{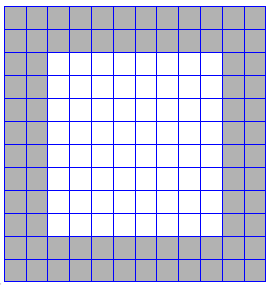}
\hspace*{8mm}
\includegraphics[height=28mm]{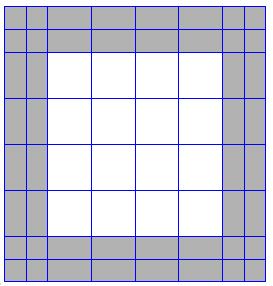}
\end{center}
\caption{Schematic display of mesh coarsening for multigrid in
  presence of PML layers in 2-D}
\label{fig:coarsening_idea}
\end{figure}

A potential solution to the second problem is to avoid mesh coarsening
and refinement in the direction normal to the PML layer, i.e.\ the
direction of the rapid variation of the coefficients $\sigma_x$, and
$\sigma_y$. This provides a simple way to avoid certain interpolation
and discretization errors in these direction of rapid
variations. Numbering the mesh cells with half-integers, assuming
$w_{\rm pml}$ cells in the PML layer. The idea is that there is no
coarsening inside the PML layers, i.e.\ for axis $j$, the cells $1/2,
\ldots, w_{\rm pml} - 1/2$ and $N_x - w_{\rm pml} +1/2 , \ldots, N_x -
1/2$ are not coarsened while the $N_x - 2 w_{\rm pml}$ interior cells
undergo standard coarsening (and similar in the $y$-direction and
$z$-direction), see Figure~\ref{fig:coarsening_idea}.

The changes to the multigrid method concern the prolongation and
restriction operators, and the coarse level discretization. We propose
to determine both in a finite element context.  

The choice of the coarse level discretization is described in detail
in the appendix.  It is such the phase speed differences with the fine
level discretization are minimized like, the discretization discussed
in \cite{StolkEtAl2014} and it is a compact stencil discretization
like required for the domain decomposition as presented here.

The prolongation and restriction operators can be written as tensor
products of one dimensional prolongation and restriction operators,
obtained by using tent finite elements. Let $i_{\rm FP}$ be the
function maps a coarse point index to the corresponding fine point
index along one of the axes, and that the function $r_{\rm C}(i)$
evaluates to ``true'' when cell $i$ is refined and ``false''
otherwise. Letting $i$ refer to any coarse mesh point and $\tilde{i} =
i_{\rm FP}(i)$ to the corresponding fine mesh point, the 1-D
prolongation operator is given by
\begin{equation}
  (P u)_{\tilde{i}} = u_i
\end{equation}
and
\begin{equation}
  (P u)_{\tilde{i} + 1} = \frac{1}{2} (u_{i} + u_{i+1})
\qquad
\text{if $r_{\rm C}(i+1/2)$.}
\end{equation}
This defines the prolongation operator.  The restriction operator is
its transpose. This concludes the description of the modified two-grid
method. 

\section{Implementation}
\label{sec:implementation}

We have developed a parallel implementation of the above described
method in three dimensions on a distributed memory machine (Linux
cluster) using MPI. The parallel implementation is fairly
straightforward, except for the domain decomposition
preconditioner. Inside the two-grid method, a Cartesian distribution
of the degrees of freedom over the compute nodes is used. 
The $\omega$-Jacobi smoother, and the restriction and prolongation
operators were implemented in a matrix-free fashion. Each time one of
these operators is applied, some communication is done between 
nodes that are neighbors in the Cartesian compute grid. 

A 2-D Cartesian compute grid is used for easy combination with the
sweeping preconditioner. Degrees of freedom are not distributed over
the sweeping axis. 

The main difficulty in the sweeping preconditioner concerns the
subdomain solves. These are done using a sparse direct solver. In the
UD-sweep all subdomains solves are done consecutively. In the X-sweep
several two solves can be done simultaneously, while in the NX-sweep
multiple solves can be done simultaneously. In particular the UD-sweep
leads to a challenging parallellization problem. 

There are several
software packages avaible to perform sparse direct solves, which allow
for various degrees of parallellization. We investigated two strategies
\begin{enumerate}
\item
Our first strategy was to use all the available compute nodes for each
solve using the Clique parallel solver of
\cite{PoulsonEtAl2013}. This solver is designed for use on many-core
systems. However, we found that solutions were sometimes incorrect. We
attribute this to limitations in the strategies for choosing pivots
(pivots were chosen inside previously chosen nested-dissection nodes.)
When these experiments were done, this solver was still in development
and the problem could be absent in later versions, but we have not
tested this. 
\item
%%% some Changes !!!
Our second strategy was to apply the method to multiple, say $n_{\rm
  RHS}$, right hand sides at the same time, and to apply the domain
decomposition preconditioner in a pipelined fashion. In the domain
decomposition step, the total number
of computational processes was divided in $n_{\rm RHS}$ groups (for
the UD-sweep) or $2n_{\rm RHS}$ groups (for the X-sweep) and each
group was responsible for a number of subdomain solves. By suitably
assigning the subdomain solves to the groups of processes, all groups
of processors could be busy at the same time (starting from step
$n_{\rm RHS}$ in the domain decomposition, when the pipeline was filled).
The factorizations and solves were done using the MUMPS parallel solver \cite{MUMPS:1}, 
version 4.10.0. For this solver it is known that it performs best when
the number of process is not too large compared to the size of the system.
A disadvantage of this method is that it leads to large memory
requirements, because of the storage required by GMRES. We
experimented with values of $n_{\rm RHS} \le 8$, at which value the
memory used for GMRES and for the subdomain factorizations were of
roughly the same size.  The outer iterative method and the two-grid
method were applied to $n_{\rm RHS}$ vectors simultaneously.
\end{enumerate}
Because of the incorrect solves in the first strategy, results will
only be given for the second strategy.

\section{Numerical experiments}
\label{sec:numerical_results}
 
In this section we study the numerical performance of the two-grid
sweeping preconditioner. The 2-D case is the easiest to study and vary
the various parameters.  We have studied problems of sizes up to $2048
\times 2048$ (for a square domain) and $4600 \times 1500$ (for the
Marmousi problem) on a laptop with 8GB memory using a Matlab
implementation.
For the three-dimensional example the parallel implementation that was
described in the previous section was used and the emphasis is on the
actual computation times.

%%% some Changes !!! 
In the numerical experiments below, the value $w_{\rm pml}$ refers to
the width of the PML layers introduced in the domain
decomposition. At the outer boundaries of the domain, sponge or PML
boundary layers are used as indicated.

\subsection{2-D experiments}

\begin{figure}
\begin{center}
\includegraphics[width=10cm]{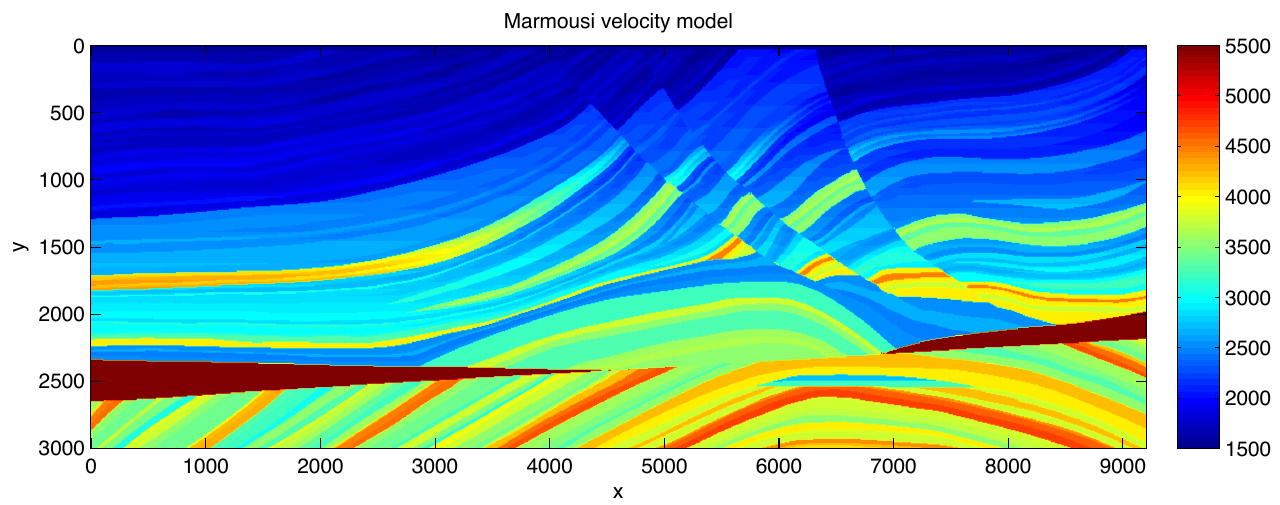}
\end{center}
\vspace*{-3ex}
\caption{Marmousi velocity model}
\label{fig:Marmousi_velocity}
\end{figure}

The first of our 2-D experiments concerns a comparison of the new
transmission conditions to those of \cite{Stolk2013} and of the new
X-sweep method with the UD sweep method used in \cite{Stolk2013}. The
comparison is done for two different discretizations, for different
values of $w_{\rm pml}$ and for two velocity models: a constant model
of size $1024 \times 1024$ grid point and the Marmousi model of size
$2300 \times 750$.  The latter model is displayed in
Figure~\ref{fig:Marmousi_velocity}. In both models a minimum of 10
points per wave length is used. Sponge boundary layers of
thickness 36 were used.  Iteration numbers to reduce the
residual by a factor $10^{-6}$ are given in
Table~\ref{tab:compare_transmission_UD_X}.

The new transmission conditions are consistent with arbitrary 9 point
discretizations, not only the standard 5 point discretization and
indeed this shows from the results. In the old transmission method,
the planar transmission source radiates not only in the direction of
the sweep, but also backward, into the added PML layer, while this is
not the case in the new method. This fact explains that for small
$w_{\rm pml}$ the new method performs better, in both discretizations.

\begin{table}
\begin{center}
\begin{tabular}{|c|c||cccc||cccc|} \hline
\multicolumn{10}{|c|}{Constant medium $1024 \times 1024$} \\ \hline
  & $n_{\rm dom}$
  & \multicolumn{4}{c||}{standard 5pt discretization} 
  & \multicolumn{4}{c|}{opt 9pt discretization} \\ \hline
  &
  & \multicolumn{2}{c}{$\;$UD-sweep$\;$}
  & \multicolumn{2}{c||}{X-sweep}
  & \multicolumn{2}{c}{$\;$UD-sweep$\;$}
  & \multicolumn{2}{c|}{X-sweep} \\
  & 
  & T1 & T2 
  & T1 & T2 
  & T1 & T2 
  & T1 & T2 \\
$w_{\rm pml} = 3$ & 78 
  &  7 & 13 &  8 & 14 &  8 & 63 &  9 & 59 \\
$w_{\rm pml} = 4$ & 60  
  &  6 &  6 &  7 &  7 &  6 & 14 &  6 & 14 \\
$w_{\rm pml} = 5$ & 49  
  &  5 &  5 &  6 &  6 &  4 & 10 &  5 & 10 \\ \hline
\multicolumn{10}{|c|}{Marmousi $2300 \times 750$} \\ \hline
$w_{\rm pml} = 3$  & 169
  & $\;$18$\;$ & $\;$58$\;$ & $\;$18$\;$ & $\;$53$\;$ 
  & $\;$18$\;$ & $\;$30$\;$ & $\;$19$\;$ & $\;$30$\;$ \\
$w_{\rm pml} = 4$  & 131
  & 12 & 12 & 12 & 14 & 11 & 25 & 12 & 26 \\
$w_{\rm pml} = 5$  & 107
  &  9 &  9 & 10 & 11 &  9 & 12 & 10 & 13 \\ \hline
\end{tabular}
\end{center}
\caption{Iteration counts for different transmission conditions for
  the UD and X-sweep preconditioners. T1 refers to the new
  transmission conditions, T2 to those of \cite{Stolk2013}.}
\label{tab:compare_transmission_UD_X}
\end{table}

We next study the two-grid method and the hybrid two-grid domain
decomposition preconditioner. To choose the smoother parameters, we
study the convergence of the two grid method with exact coarse level
inverse. We vary $\nu$ (the number of pre- and postsmoothing steps)
and $\omega_{\rm Jac}$, the relaxation constant. The model is the unit
square with unit velocity discretized with $1536 \times 1536$ points
(excluding sponge or PML layers) and with frequency
$\frac{\omega}{2\pi} = 153.6$ (10 points per wavelength). This is
about the largest problem that can still be done without using
excessive amounts of swap memory. The tests are done using sponge 
boundary layers of thickness 36 and PML layers of thickness 4. The
results in Table~\ref{tab:iter_smoothing_parameters} show that
$\nu = 3$ and $\omega_{\rm Jac} = 0.8$ gives good results. The 
improvements in iteration count found for even larger values of $\nu$ are not
found in other experiments involving domain decomposition.
Therefore we choose $\nu = 3$ and $\omega_{\rm Jac} = 0.8$ for the 2-D
problem. Good results are obtained for both sponge and PML layers, we
will study the difference further in other examples.
\newlength{\myts} \setlength{\myts}{0.7ex}
\begin{table}
\newcommand{\gtrinfo}{\multicolumn{2}{c}{$>99$}}
\newcommand{\gtrinfolast}{\multicolumn{2}{c|}{$>99$}}
\begin{center}
\begin{tabular}{|c|r@{\hskip \myts}lr@{\hskip \myts}l%
r@{\hskip \myts}lr@{\hskip \myts}lr@{\hskip \myts}l|} \hline
\multicolumn{11}{|c|}{Sponge} \\
\hline
    & \multicolumn{2}{c}{$\omega{\rm Jac} = 0.5$}
    & \multicolumn{2}{c}{0.6}
    & \multicolumn{2}{c}{0.7}
    & \multicolumn{2}{c}{0.8}
    & \multicolumn{2}{c|}{0.9} \\ \hline
$\nu = 1$ 
    & \gtrinfo   & 69&(114) & 47&(77) & 26&(43) & 26&(43) \\
  2 &  $\;\;$28&(54) & 18&(36)  & 13&(26) &  9&(19) & 12&(25) \\
  3 &    15&(34) &  9&(21)  &  7&(17) &  5&(12) &  9&(21) \\
  4 &     9&(24) &  7&(19)  &  5&(14) &  5&(14) &  8&(22) \\
  5 &     7&(21) &  5&(16)  &  4&(13) &  4&(13) &  7&(21) \\
  6 &     6&(20) &  5&(18)  &  4&(15) &  4&(15) &  6&(20) \\
\hline
\multicolumn{11}{|c|}{PML} \\
\hline
    & \multicolumn{2}{c}{$\omega{\rm Jac} = 0.5$}
    & \multicolumn{2}{c}{0.6}
    & \multicolumn{2}{c}{0.7}
    & \multicolumn{2}{c}{0.8}
    & \multicolumn{2}{c|}{0.9} \\ \hline
$\nu = 1$ 
    & \gtrinfo    & \gtrinfo   & \gtrinfo   & 54&(85) & \gtrinfolast \\
  2 &    67&(120) &    30&(55) &    18&(33) & 15&(28) & \gtrinfolast \\
  3 &    23&(48)  &    14&(30) &    10&(22) & 10&(23) &    80&(161) \\
  4 &    13&(32)  &    10&(25) &     9&(22) & 10&(25) &    70&(158) \\
  5 &    10&(28)  &     9&(25) &     9&(25) & 10&(28) &    70&(175) \\
  6 &     9&(28)  &     9&(27) &    10&(30) & 12&(37) &    82&(224) \\
\hline
\end{tabular}
\end{center}
\caption{Number of iterations (computation times for the solve phase) 
as a function smoothing parameters for a constant velocity, 10 points
per wavelength, and mesh size 1536$\times$1536.}
\label{tab:iter_smoothing_parameters}
\end{table}

Next we study the convergence for different values of $w_{\rm pml}$
and the problem size. We also include the exact coarse scale
solver. This is done for two problems, the constant-velocity unit
square and the Marmousi model. For the constant velocity model,
10 points per wavelength fine scale discretization was used.
The values of $S_{\rm pml}$ are chosen to be 15, 20 and 25
respectively for $w_{\rm pml} = 3$, $4$ and $5$. 
For the outer boundaries sponge boundary layers of thickness 36
and PML layers of thickness 4 were used. We determined
iteration counts and the time for the solve phase. Setup times were 
of the same order of magnitude as the solve times. 
Results are in Table~\ref{tab:iter_constsquare_and_Marmousi}.
The number of subdomains used, given by $\left\lfloor
\frac{N_x}{2w_{\rm pml}+1} \right\rfloor$, depended on $w_{\rm pml}$ 
and is also indicated in the table (in the column labeled $n_{\rm dom}$).
It is clearly seen that for larger problems also a larger value of 
$w_{\rm pml}$ should be used because the number of iterations grows
faster than the extra cost of thicker PML layers.
In some examples good convergence was obtained  
using up to 250 subdomains.

Next we test the X-sweep, and the NX-sweep approaches described in
section~\ref{sec:dsdd_method}, involving simultaneous and partial
sweeps. Iteration numbers for these approaches for our largest
constant and Marmousi examples are given in
Table~\ref{tab:iter_sweep_patterns}. In both cases we see that the
UD-sweep pattern can be replaced by the X-sweep pattern at little or
no cost. The method with partial sweeps performs poorly. The gain in
computation time that can be obtained by performing the partial sweeps
in parallel disappears because of the additionally
required iterations.

\begin{table}
\begin{center}
\begin{tabular}{|ccc|r@{\hskip 0.7ex}lr@{\hskip 0.7ex}lr@{\hskip 0.7ex}lr@{\hskip 0.7ex}l|} \hline
\multicolumn{11}{|c|}{Constant medium, PML} \\ \hline
size 
  & freq.\
  & $n_{\rm dom}$ for
  & \multicolumn{2}{c}{exact} 
  & \multicolumn{2}{c}{$w_{\rm pml} = 3$} 
  & \multicolumn{2}{c}{$w_{\rm pml} = 4$} 
  & \multicolumn{2}{c|}{$w_{\rm pml} = 5$} 
\\
  &
  & $w_{\rm pml}=$3/4/5
  &&
  &&
  &&
  &&
\\ \hline
$ 256 \times  256$ & 25.6  & 23/18/14
  & 10&(0.51)& 10&(0.75)& 10&(0.75)& 10&(0.76)\\
$ 512 \times  512$ & 51.2  & 41/32/26
  &  10&(2.4) & 11&(3.3) & 10&(2.9) & 10&(2.9)\\
$1024 \times 1024$ & 102.4 & 78/60/49
  & 10&(10)  & 13&(14)  & 11&(13)  & 10&(12) \\
$2048 \times 2048$ & 204.8 & 151/117/96
  & 11&(61)  & 32&(129) & 15&(64)  & 12&(55) \\ \hline
\multicolumn{11}{|c|}{Constant medium, sponge} \\ \hline
$ 256 \times  256$ &  25.6  & 23/18/14
  & 5&(0.46) &  5&(0.63) & 5&(0.63) & 5&(0.70)\\
$ 512 \times  512$ &  51.2  & 41/32/26 
  & 5&(1.6)  &  6&(2.3)  & 5&(1.9)  & 5&(1.9)\\
$1024 \times 1024$ & 102.4 & 78/60/49 
  & 5&(5.8)  &  7&(8.8)  & 6&(7.9)  & 6&(7.9)\\
$2048 \times 2048$ & 204.8 & 151/117/96 
  & 6&(41)   & 10&(46)   & 7&(38)   & 7&(34) \\ \hline
\multicolumn{11}{|c|}{Marmousi model, PML} \\ \hline
$ 575 \times  188$ & 9.4  & 46/36/29
  & 13&(1.2) & 14&(1.8) & 14&(1.8) & 14&(1.8) \\
$1150 \times  375$ & 18.8 & 87/67/55
  & 15&(5.8) & 14&(6.6) & 14&(6.6) & 14&(6.6) \\
$2300 \times  750$ & 37.5 & 169/131/107
  & 13&(20)  & 17&(29)  & 14&(25)  & 14&(25)  \\
$4600 \times 1500$ &  75  & 333/259/212
  & 12&(*)   & 39&(*)   & 17&(*)   & 14&(*)   \\ \hline
\multicolumn{11}{|c|}{Marmousi model, sponge} \\ \hline
$ 575 \times  188$ & 9.4  & 46/36/29
  & 10&(1.4) & 10&(1.9) & 10&(1.9) & 10&(1.9) \\
$1150 \times  375$ & 18.8 & 87/67/55
  & 12&(5.9) & 13&(7.8) & 13&(8.1) & 13&(7.9) \\
$2300 \times  750$ & 37.5 & 169/131/107
  & 11&(22)  & 14&(29)  & 13&(26)  & 13&(26)  \\
$4600 \times 1500$ &  75  & 333/259/212
  & 10&(*)   & 25&(*)   & 14&(*)   & 13&(*)   \\ \hline
\end{tabular}
\end{center}
\caption{Iteration numbers (and time per solve in seconds) as a 
  function of $w_{\rm pml}$ and problem size for a constant velocity
  and the Marmousi model. (*)
  denotes long times, between 200 and 600 seconds, due to shortage of RAM}
\label{tab:iter_constsquare_and_Marmousi}
\end{table}

\begin{table}
\begin{center}
\begin{tabular}{|c|cc|} \hline
velocity    & CONSTANT & MARMOUSI \\
size        & $2048 \times 2048$ 
                       & $4600 \times 1500$ \\ \hline
UD-sweep    &  8       & 15 \\
X-sweep     &  8       & 16 \\
NX-sweep(2) & 40       & 54 \\
NX-sweep(4) & 46       & 68 \\
NX-sweep(8) & 63       & 99 \\ \hline
\end{tabular}
\end{center}
\caption{Iteration numbers as a function of sweep type for the
  constant and Marmousi velocity models. For the NX-sweep pattern the
  number $N_{\rm cell}$ is indicated between the brackets.}
\label{tab:iter_sweep_patterns}
\end{table}

\subsection{The 3-D SEG-EAGE salt model}

\begin{figure}
\begin{center}
(a) \hspace*{6cm} (b)\\
\includegraphics[width=6cm,height=29mm]{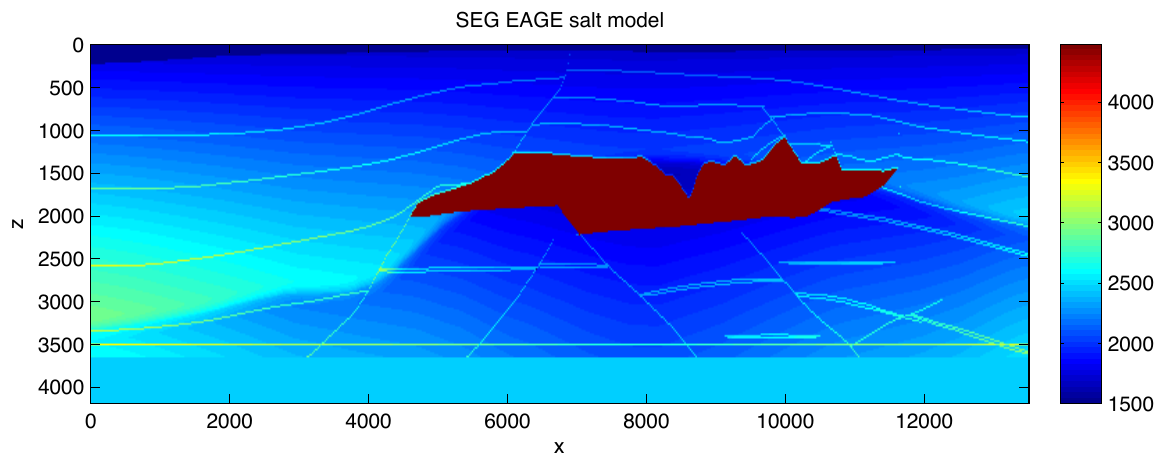}
\includegraphics[width=6cm,height=29mm]{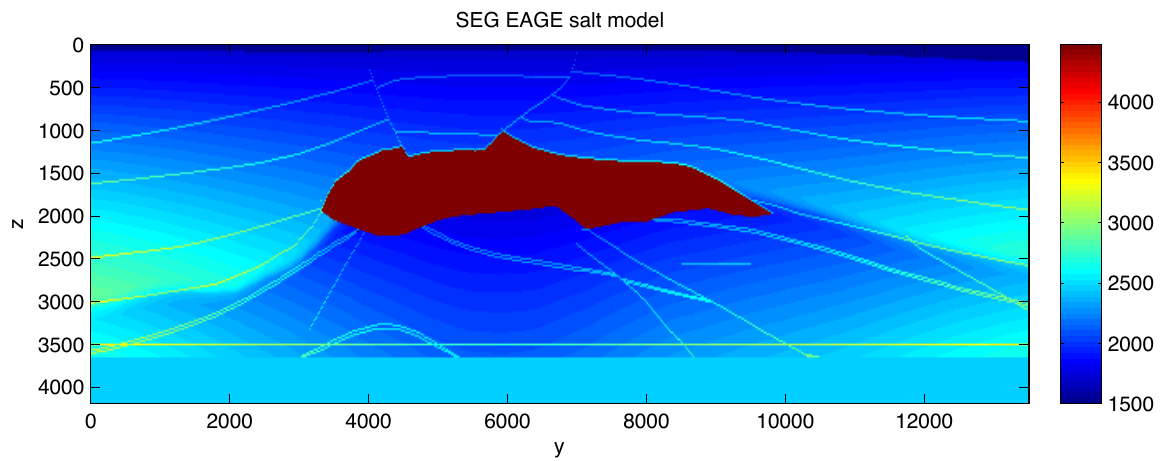}

\smallskip

(c) \hspace*{6cm} (d)\\
\includegraphics[width=6cm,height=29mm]{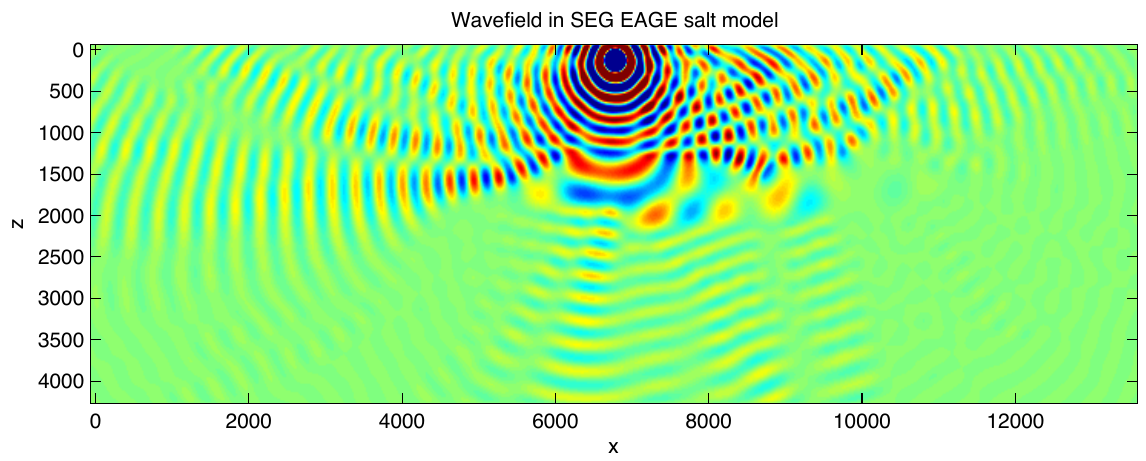}
\includegraphics[width=6cm,height=29mm]{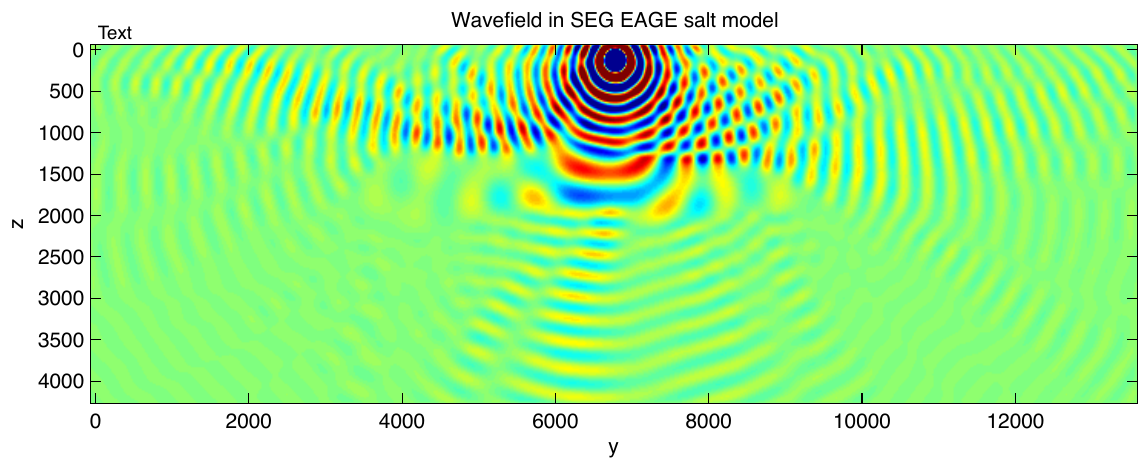}
\end{center}
\caption{SEG-EAGE salt velocity model: (a) $(x,z)$ slice at $y=6740$ 
(b) $(y,z)$ slice at $x=6740$. Solution to the Helmholtz equation at 
$7.5$ Hz: (c) $(x,z)$ slice at $y=6740$ 
(d) $(y,z)$ slice at $x=6740$.}
\label{fig:saltmodel_pictures}
\end{figure}

The SEG-EAGE salt model is a 3-D synthetic Earth model from exploration
geophysics. The original model is of size 13500 x 13500 x 4200 meter,
discretized with 20 m grid spacing. We apply the two-grid sweeping
preconditioner to solve the Helmholtz equation with this velocity at
four different frequencies from $3.75$ to $7.5$ Hz, using a minimum of 
10 points per wave length. At the outer boundaries, we used PML boundary layers
of width 3 grid points. In the domain decomposition, 
we used $w_{\rm pml} = 3$. Three
iterations of $\omega$-Jacobi with $\omega = 0.6$ were used as
smoother in the two-grid method. The
right hand side was chosen randomly. Convergence for the random right hand side
typically required about 1 iteration extra compared to the point source.
Slices of the model, and a solution with a points source of the
Helmholtz equation at 7.5
Hz are displayed in Figure~\ref{fig:saltmodel_pictures}.
The problem studied has about $1.0 \cdot 10^8$ degrees of freedom.

Computations were done the Lisa cluster at surfsara (www.surfsara.nl)
using the implementation described in
section~\ref{sec:implementation}. For parallel computations this
systems contains 32 nodes with each two intel Xeon processors 
E5-2650 v2 running at 2.60 GHz and 64 GB memory, connected by Mellanox 
FDR Infiniband. The use of two intel Xeon units results in 16 cores per node. 
A maximum of 16 nodes were used in parallel for these computations.

As described in section~\ref{sec:implementation}, the algorithm solves
multiple right hand sides at the same time, using subgroups of
processes for the subdomain solves in combination with pipelining. The
number of right hand sides was chosen $\le 8$, to control the memory
use. The size of the subgroups was varied between 8 and 32. For larger
subgroups, larger problems can be solved using the parallel algorithm.

Results, in particular iteration counts and computation times, of the
computations are given in Table~\ref{tab:SaltModel}.  Our main
conclusion is that there is large improvement in computation times
and memory use compared to the pure sweeping methods described in
\cite{PoulsonEtAl2013}, such that the method becomes comparable to in
computation times to some of the fastests methods in the literature,
see for example \cite{CalandraEtAl2013}, where a combination of a
two-grid and a shifted Laplacian method was considered and
\cite{RiyantiEtAl2007,WangDeHoopXia2011} for further 
examples of solvers applied to large scale examples.

Considering the results as a function of problem size we see that
computation times increase with problem size, even if the number of
processes also increases. Several factors contribute to this: the
number of iterations increases slowly, the cost of the sparse direct
solve increases somewhat faster than linearly and cost related to the
parallellization will also typically increase. When the MUMPS solver
is used with 32 cores, the computation times are somewhat longer
compared to 8 or 16 cores. While it is difficult to explain this
precisely, it is likely that the slow communication over
multiple nodes (instead of just within a node) contributes to this.

\begin{table}
\begin{center}\footnotesize
\begin{tabular}{|l|c|c|c|c|}
\hline
Freq.\ (Hz)
  & 3.75
  & 4.72
  & 5.95
  & 7.5
\\
Problem size 
  & 338x338x106
  & 426x426x132
  & 536x536x166
  & 676x676x210
\\
\#layers
  & 25
  & 30
  & 40
  & 48
\\
\#dof
  & $1.3 \cdot 10^7$
  & $2.5 \cdot 10^7$
  & $5.0 \cdot 10^7$
  & $1.0 \cdot 10^8$
\\
Cores 
  & 32
  & 64
  & 128
  & 256
\\ \hline\hline
\multicolumn{5}{|c|}{UD-SWEEP} 
\\ \hline 
iterations
  & 11
  & 12
  & 12
  & 14
\\ \hline
\multicolumn{5}{|c|}{Mumps 16 cores} 
\\ \hline
\#rhs
  & 2
  & 4 
  & 8
  & 
\\
setup time (s) 
  & 47
  & 54
  & 66
  &
\\
solvetime/rhs
  & 27
  & 26
  & 44
  &
\\ \hline 
\multicolumn{5}{|c|}{Mumps 32 cores} 
\\ \hline
\#rhs
  & 1
  & 2 
  & 4
  & 8
\\
setup time (s) 
  & 74
  & 82
  & 94
  & 144
\\
solvetime/rhs
  & 36
  & 48
  & 52
  & 67
\\ \hline\hline
\multicolumn{5}{|c|}{X-SWEEP} 
\\ \hline 
iterations
  & 11
  & 12
  & 13
  & 15
\\ \hline
\multicolumn{5}{|c|}{Mumps 8 cores} 
\\ \hline
\#rhs
  & 2
  & 4 
  & 8
  & 
\\
setup time (s) 
  & 39
  & 44
  & 62
  &
\\
solvetime/rhs
  & 20
  & 26
  & 39
  &
\\ \hline
\multicolumn{5}{|c|}{Mumps 16 cores} 
\\ \hline
\#rhs
  & 1
  & 2 
  & 4
  & 8
\\
setup time (s) 
  & 49
  & 54
  & 66
  & 96
\\
solvetime/rhs
  & 22
  & 27
  & 31
  & 62
\\ \hline
\multicolumn{5}{|c|}{Mumps 32 cores} 
\\ \hline
\#rhs
  & 
  & 1 
  & 2
  & 4
\\
setup time (s) 
  &
  & 82
  & 86
  & 107
\\  
solvetime/rhs
  &
  & 43
  & 60
  & 80
\\ \hline
\end{tabular}
\end{center}
\caption{Simulation results for the 3-D SEG-EAGE salt model}
\label{tab:SaltModel}
\end{table}

\section{Discussion}
\label{sec:discussion}

In this work we used a two-grid method to accellerate a Helmholtz
solver based on a sweeping preconditioner. This resulted in a new
method that we call two-grid sweeping preconditioner. A priori it was
not clear that such a method would work, as both the sweeping
preconditioner and the two-grid method are used in new conditions.

With the two-grid method as outer method, the cost of the sweeping
preconditioner is strongly reduced. When problems of the same size are
considered, computation times appear to be roughly comparable to those
of the method of \cite{CalandraEtAl2013}, where a combination of a
two-grid and a shifted Laplacian method were considered. Thus the
methods is comparable in performance to some of the fastests methods
in the literature.
(See \cite{PoulsonEtAl2013,RiyantiEtAl2007,WangDeHoopXia2011} 
for other works that consider large scale examples.)

Parallellization of the numerical linear algebra remains a challenge
for these methods. The performance of sweeping preconditioners is
determined in part by the possibilities and limitations of parallel
solvers like MUMPS \cite{MUMPS:1} and Clique
\cite{PoulsonEtAl2013}. For reasons explained in
section~\ref{sec:implementation} we used MUMPS. The version which was
used doesn't scale very well to large numbers of processes.
Improvements in this area will be useful for large scale parallel
application of the methods.

If $w_{\rm pml}$ and the thickness of the layers is kept fixed, the
preconditioner can be applied with cost log-linear in the number of
unknowns, because for a single layer of size $n \times n \times d$,
the cost for solving the factorized system is $O(d^2 n^2 \log n)$, 
cf.\ \cite{EngquistYing2011,George1973}).  The numerical
results show that quite small values of $w_{\rm pml}$ can be used
(e.g.\ $w_{\rm pml} = 3$ with more than 100 subdomains).  However, we
find that to keep good convergence for larger number of subdomains, 
$w_{\rm pml}$ should increase slowly with problem size.

\bibliographystyle{abbrv}
\bibliography{helmmgsp}

\appendix

\section{A Helmholtz discretization for use in the two-grid sweeping
  preconditioner with PML boundary layers}

In this section we discuss a discretization that can be used on meshes
of the type displayed in Figure~\ref{fig:coarsening_idea}, where
inside the PML layers, the coarsening only takes place in the
tangential directions.  This is
done using a variant of the multigrid finite element method. The
result can be used as coarse level discretization in a multigrid
method, as explain in section~\ref{sec:two_grid_pml}. The construction
of a coarse level operator with phase speeds matching those of the
fine level operator is achieved using the equivalence between finite
element schemes with general testfunctions and finite
difference schemes. This allows us to reproduce the behavior of the
optimized finite difference method of \cite{StolkEtAl2014} in the
current setting.  We will treat the 3-D case, which is slightly more
complicated than the the 2-D case.

The discretization is done using rectilinear (product) meshes
with mesh points $(x_i,y_j,z_k)$, $0 \le i \le N_x, 0 \le j \le N_y$
and $0 \le k \le N_z$. The cells will be numbered such that cell
$i+1/2$ is between points $i$ and $i+1$. Cell size parameters of cell
$(i+1/2,j+1/2,k+1/2)$ are $h_{1,i+1/2}$, $h_{2,j+1/2}$ and
$h_{3,k+1/2}$. This allows for regular and non-regular meshes A
regular mesh of this type can be used for finite differences.  For a
regular mesh, $h$ will denote the mesh parameter. General rectilinear
meshes of this type can be used for finite element discretizations. We
assume the nodes are the eight corners of each cell, and degrees of
freedom are denoted by $u_{i,j,k}$. The degrees of freedom are located
at points with $1 \le i \le N_x -1$, $1 \le j \le N_y -1$, $1 \le k
\le N_z -1$ because Dirichlet boundary conditions are used.

In the remainder of this section we first revisit the optimized finite
differences from \cite{StolkEtAl2014}. We then describe a general
finite element discretization. In the third subsection we describe how
to choose coefficients in this general finite element discretization
to recover the optimized finite differences in the regular, non-PML
part of the mesh. This yields the discretization that we used in the
two-grid method when PML layers were present. In the last subsection
of this appendix we present a further result on the connection between
finite elements and optimized finite differences.

\subsection{Optimized finite differences}
\label{sec:optimized_fd}
Optimized finite differences for frequency domain simulation in the
plane are described for example in \cite{JoShinSuh1996}. In
\cite{StolkEtAl2014} a different version was introduced for both two
and three dimensions which was applied in a multigrid method.  See
also \cite{BabuskaEtAl1995} and further discussion
in \cite{Stolk2016_Dispersion_minimizing_scheme}.  We will explain in detail the 3-D
method of \cite{StolkEtAl2014}, the 2-D version is derived in the same way.

First we define some discrete operators. Define $M_j$, $j=0,1,2,3$ by
\begin{equation}
\begin{split}
  ( M_0 u )_{i,j,k} = {}& u_{i,j,k} 
\\
  ( M_1 u )_{i,j,k} = {}& \frac{1}{6} \big(
u_{i-1,j,k}+u_{i+1,j,k}+u_{i,j-1,k}+u_{i,j+1,k}+u_{i,j,k-1}+u_{i,j,k+1} \big)
\\
  ( M_2 u )_{i,j,k} = {}& \frac{1}{12} \big(
  u_{i-1,j-1,k}+u_{i+1,j-1,k}+u_{i-1,j+1,k}+u_{i+1,j+1,k}
  +u_{i-1,j,k-1}+u_{i+1,j,k-1}
\\
  {}& +u_{i-1,j,k+1}+u_{i+1,j,k+1}
  +u_{i,j-1,k-1}+u_{i,j+1,k-1}+u_{i,j-1,k+1}+u_{i,j+1,k+1} \big)
\\
  ( M_3 u )_{i,j,k} = {}& \frac{1}{8} \big(
  u_{i-1,j-1,k-1}+u_{i+1,j-1,k-1}+u_{i-1,j+1,k-1}+u_{i+1,j+1,k-1}
\\
  {}& +u_{i-1,j-1,k+1}+u_{i+1,j-1,k+1}+u_{i-1,j+1,k+1}+u_{i+1,j+1,k+1} \big).
\end{split}
\end{equation}
All of these are second order discretizations of the identity operator
in 3-D. Similarly, for 2-D field $u_{i,j}$, consider the operators
$N_j$, $j=0,1,2$ given by
\begin{equation}
\begin{split}
  ( N_0 u )_{i,j} = {}& u_{i,j} 
\\
  ( N_1 u )_{i,j} = {}& 
\frac{1}{4} \big( u_{i-1,j}+u_{i+1,j}+u_{i,j-1}+u_{i,j+1} \big)
\\
  ( N_2 u )_{i,j} = {}& 
\frac{1}{4} \big( u_{i-1,j-1}+u_{i+1,j-1}+u_{i-1,j+1}+u_{i+1,j+1} \big)
\end{split}
\end{equation}
These form discretizations of the identity operator in 2-D. By
$N_a^{(l,m)}$ we denote these operators acting along the $(x_l,x_m)$
axes. Furthermore, denote by $D_{2,{\rm FD}}$ the discrete second
order derivative
\begin{equation}
  (D_{2,{\rm FD}} u)_i = \frac{1}{h^2} \left( u_{i-1} - 2u_i + u_{i+1} \right) .
\end{equation}
By $D_2^{(l)}$ we denote this operator acting along the $x_l$ axis.

We will next define a five parameter family of second order discrete 
Helmholtz operators. Given 5 coefficients $c_j$, $j = 1,2,3,4,5$, denote
\begin{equation}
\begin{split}
  \widetilde{M} = {}& c_1 M_0 u+ c_2 M_1 u + c_3 M_2 u + (1-c_1-c_2-c_3) M_3 u
\\
  \widetilde{N} = {}& c_4 N_0 + c_5 N_1 + (1-c_4-c_5) N_2 .
\end{split}
\end{equation}
By $\widetilde{N}^{(l,m)}$ we will denote versions of these operator
acting along the $(x_l,x_m)$ axes.
The operators $\widetilde{M}$, $\widetilde{N}$ are weighted average of
second order discretizations of the identity, and are hence second
order discretizations of the identity themselves.
We use them to define a five parameter family of second order 
discretizations of the Helmholtz operator, with a compact 
$3 \times 3 \times 3$ stencil as follows
\begin{equation} \label{eq:define_H_FDOPT}
\begin{aligned}
  (H_{\rm FDOpt} u)_{i,j,k}
  \stackrel{\rm def}{=}
{}& - k_{i,j,k}^2 (\widetilde{M} u)_{i,j,k}
  - ( (D_{2,{\rm FD}}^{(1)} \otimes \widetilde{N}^{(2,3)}) u)_{i,j,k}  
  - ( (D_{2,{\rm FD}}^{(2)} \otimes \widetilde{N}^{(1,3)}) u)_{i,j,k}  
\\
  {}& - ( (D_{2,{\rm FD}}^{(3)} \otimes \widetilde{N}^{(1,2)}) u)_{i,j,k}  
\\
   = {}&f_{i,j,k} .
\end{aligned}
\end{equation}
In 2-D, a similar formula can be made with three independent
coefficients $c_j$, $j=1,2,3$. 

We now have five coefficients that can be chosen (or three in 2-D).
The phase speed of the numerical method depends on the product $k h$,
or equivalently on the number of points per wavelength $G = \frac{2
  \pi}{h k}$ and on the direction of the wave that is considered. In
addition it depends on the choice of the coefficients $c_l$. In
\cite{JoShinSuh1996} the coefficients $c_l$, $l=1,2,3$ for the 2-D
case, were fixed so as to minimize the maximum of the absolute 
difference between the exact and the numerical phase speeds
(to be precise, Jo Shin Suh considered a different set of basic 
operators and an equivalent set of coefficients was fixed). Here 
the maximum was taken over all angles and $G \ge 4$. In this way, 
a numerical method with much better dispersion properties than 
standard second order finite differences was obtained.

\begin{figure}
\begin{center}
(a)\\
\includegraphics[width=80mm]{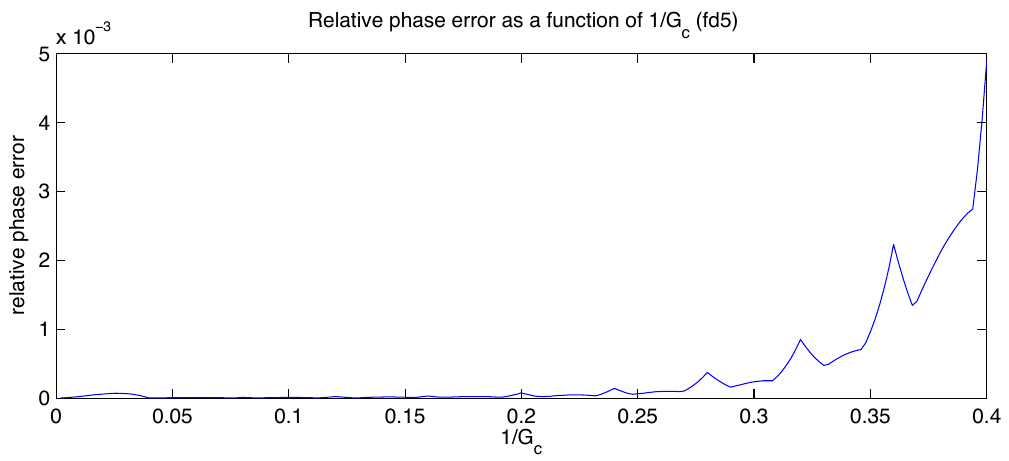}\\
(b)\\
\includegraphics[width=80mm]{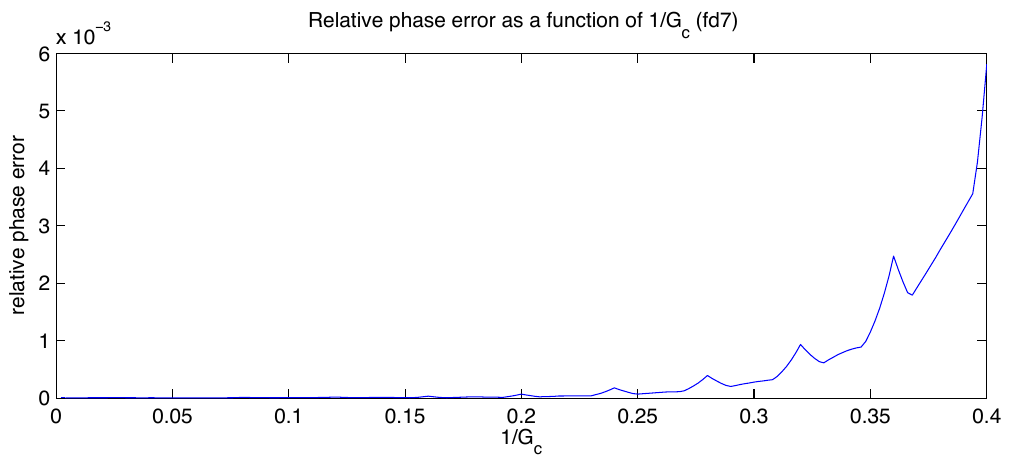}
\end{center}
\caption{Fine-coarse phase speed error using optimized finite
  differences, maximum over angle (a) two dimensions; (b) three dimensions.}
\label{fig:error_optFD}
\end{figure}

Stolk et al.\ \cite{StolkEtAl2014}
observed that the phase speed errors can be further reduced if $c_j$
depends on $1/G$ (using $1/G$ is slightly more convenient than $G$).  To
represent the functions $c_j(1/G)$ simple linear interpolation was
chosen. I.e.\ the function $c_j(1/G)$ was parameterized by support
points $1/G_k$, and values $c_j(1/G_k)$, and given by linear
interpolation for values of $1/G$ between the support points. An 
optimization procedure was done to find values $c_j(1/G_k)$ such 
phase speed differences between the coarse and fine scale
methods of a two-grid method were minimal over the considered range of
$1/G$. The values $1/G_k$ and $c_l(1/G_k)$ for
both the 2-D and 3-D case are given in Table~\ref{tab:optFD_coeff}.
Graphs of the error (maximum over angle) are given in 
Figure~\ref{fig:error_optFD}. In this way the phase speed differences between the fine
and coarse scale methods could be reduced very strongly, to about $2
\cdot 10^{-4}$ for $G \ge 4$.

\begin{table}
\begin{center}
(a)\\
\begin{tabular}{|r|ccc|} \hline \small
$1/G_k$  & $c_1$ & $c_2$ & $c_3$ \\ \hline
0.00 & 0.61953 & 0.45295 & 0.77363 \\
0.04 & 0.63691 & 0.47535 & 0.87242 \\
0.08 & 0.62988 & 0.48633 & 0.86400 \\
0.12 & 0.62610 & 0.48880 & 0.84984 \\
0.16 & 0.62289 & 0.48759 & 0.83017 \\
0.20 & 0.62596 & 0.47106 & 0.80852 \\
0.24 & 0.62213 & 0.46478 & 0.78215 \\
0.28 & 0.61036 & 0.47016 & 0.74857 \\
0.32 & 0.59107 & 0.48468 & 0.70553 \\
0.36 & 0.56369 & 0.50746 & 0.65062 \\
0.40 & 0.52412 & 0.54163 & 0.57676 \\
\hline
\end{tabular}

\medskip

(b)\\
\begin{tabular}{|r|ccccc|}\hline \small
$1/G_k$  & $c_1$ & $c_2$ & $c_3$ & $c_4$ & $c_5$ \\ \hline
0.00 & 0.56428 & 0.35970 & 0.20490 & 0.77998 & 0.17505 \\
0.04 & 0.56571 & 0.36071 & 0.20541 & 0.78635 & 0.17442 \\
0.08 & 0.56298 & 0.36150 & 0.20719 & 0.78273 & 0.16881 \\
0.12 & 0.56540 & 0.35620 & 0.20287 & 0.76438 & 0.18678 \\
0.16 & 0.56370 & 0.35299 & 0.20299 & 0.74684 & 0.19603 \\
0.20 & 0.55813 & 0.35277 & 0.20452 & 0.72755 & 0.20131 \\
0.24 & 0.54673 & 0.35830 & 0.20693 & 0.70298 & 0.20847 \\
0.28 & 0.52423 & 0.38368 & 0.19633 & 0.66863 & 0.22424 \\
0.32 & 0.49946 & 0.39740 & 0.20725 & 0.62734 & 0.23845 \\
0.36 & 0.47567 & 0.40216 & 0.22132 & 0.58198 & 0.25329 \\
0.40 & 0.45011 & 0.36784 & 0.29962 & 0.53417 & 0.23589 \\
\hline
\end{tabular}
\end{center}
\caption{Coefficients for optimized finite differences 
with phase speeds matching those of standard second order 
finite differences (a) two dimensions, (b) three dimensions.}
\label{tab:optFD_coeff}
\end{table}

\subsection{A class of finite element discretizations}
\label{sec:discretization}
The weak form of the Helmholtz equation with PML boundary layers reads, using that
$u$ and $v$ vanish on the boundary,
\begin{equation} \label{eq:weak_form}
  \int_\Omega \bigg[ \sum_{j=1}^3 \frac{\alpha_j^2}{\alpha_1 \alpha_2 \alpha_3} 
    \pdpd{u}{x_j}\pdpd{v}{x_j} 
  - \frac{k^2}{\alpha_1 \alpha_2 \alpha_3} u v 
  -  \frac{1}{\alpha_1 \alpha_2 \alpha_3} f v \bigg] \, dx = 0 
\end{equation}
for all $v$, where $\alpha_j$ is as defined below (\ref{eq:define_pml_fd}).

%Terminology: trial function, test function basis functions, shape functions.
To obtain a finite element method we must describe the spaces of trial
and test functions.
The trial functions associated with the nodes of the mesh and are
derived from standard trilinear shape function. I.e.\, on the unit
cube the shape function associated with the origin is
\begin{equation}
  \psi_{0,0,0} = (1 - x_1) (1 - x_2) (1 - x_3)
\end{equation}
For the test functions we will only assume that they derive in the
usual way from a single shape function $\tilde{\psi}_{0,0,0}$ on a
reference cell that is continuous and piecewise $C^1$ and symmetric
under permutation of the axes.

We assume that $k(x)$ and the $\alpha_j$ are cellwise constant. This
implies that only a few integrals of the test and trial
functions and their derivatives need to be known. 

Next we obtain an expression for the mass matrix, i.e.\ the matrix
with elements
\begin{equation}
  M_{{\rm FE},i,j,k,\tilde{i},\tilde{j},\tilde{k}} 
  = \int \frac{k^2}{\alpha_1 \alpha_2 \alpha_3}
  u_{\tilde{i},\tilde{j},\tilde{k}} v_{i,j,k} \, dx .
\end{equation}
Define
\begin{equation} \label{eq:define_I_s1s2s3}
  I_{s_1,s_2,s_3} =\int_{[0,1]^3} \psi_{0,0,0} \tilde{\psi}_{s_1,s_2,s_3} dx .
\end{equation}
Due to the symmetries there are four independent values, namely those with
$(s_1,s_2,s_3) \in \{ (0,0,0), (1,0,0), (1,1,0), (1,1,1) \}$.
We hence set
\begin{equation} \label{eq:define_I_S}
\begin{aligned}
  I_{0} = {}& I_{0,0,0} \qquad
&
  I_{1} = {}& I_{1,0,0} 
\\
  I_{2} = {}& I_{1,1,0} 
&
  I_{3} = {}&I_{1,1,1} .
\end{aligned}
\end{equation}
To easily list the contributions to the matrix we define the index sets\\
\begin{equation}
  S(i,\tilde{i}) = \left\{ \begin{array}{ll}
 \{ -1/2 \}      & \text{if $\tilde{i} = i-1$}\\
 \{ -1/2, 1/2 \} & \text{if $\tilde{i} = i  $}\\
 \{ 1/2 \}       & \text{if $\tilde{i} = i+1$}\\
  \emptyset      & \text{otherwise.}
\end{array} \right.
\end{equation}
With these definitions, we have the following expression for the mass matrix
\begin{equation}
  M_{{\rm FE},i,j,k,\tilde{i},\tilde{j},\tilde{k}}
   = \sum_{(s_1,s_2,s_3) \in S(i,\tilde{i}) \times S(j,\tilde{j}) \times S(k,\tilde{k})}
    h_{i+s_1} h_{j+s_2} h_{k+s_3}
    I_{| \tilde{i} - i | + | \tilde{j} - j | + | \tilde{k} - k |}
    \frac{k^2_{i+s_1,j+s_2,k+s_3}}{\alpha_{1,i+s_1} \alpha_{2,j+s_2}
      \alpha_{3,k+s_3}} ,
\end{equation}
where, as usual, the sum over an empty index set is zero.  As
expected, nonzero matrix elements occur when
$\max(| \tilde{i} - i |, | \tilde{j} - j |, | \tilde{k} - k | ) \le
1$. The sum is over 8, 4, 2, or 1 cells, depending whether the vector
$(\tilde{i} - i,\tilde{j} - j,\tilde{k} - k)$ is in the center, face-center, 
edge-center or vertex position of the 27 point cube $\{-1,0,1\}^3$.

By the stiffness matrix we mean the matrix whose
$(i,j,k;\tilde{i},\tilde{j},\tilde{k})$ element is given by
\begin{equation}
  \sum_{l=1}^3 \int \frac{\alpha_l^2}{\alpha_1 \alpha_2 \alpha_3} 
    \pdpd{u_{\tilde{i},\tilde{j},\tilde{k}}}{x_l}
    \pdpd{v_{i,j,k}}{x_l} \, dx .
\end{equation}
Each summand is an integral over multiple cells, and for each summand,
and each cell, the integral can be reduce to a multiple of one of the
following integrals
\begin{equation}
  J^{(l)}_{s_1,s_2,s_3} = \int_{[0,1]^3} \pdpd{\psi_{0,0,0}}{x_l}
  \pdpd{\tilde{\psi}_{s_1,s_2,s_3}}{x_l} dx ,
\end{equation}
where the $s_j$ are 0 or 1.
Taking $l = 1$, the derivative $\pdpd{\psi_{0,0,0}}{x_1} 
= (1 - x_2) (1 - x_3)$ is independent of $x_1$ and the integral 
reduces to a sum of surface integrals
\begin{equation}
\begin{split}
  J^{(1)}_{s_1,s_2,s_3} 
  = {}& - \iint \tilde{\psi}_{s_1,s_2,s_3}(1,x_2,x_3) (1 - x_2) (1 - x_3) \, dx_2 \, dx_3
\\
  {}& + \iint \tilde{\psi}_{s_1,s_2,s_3}(0,x_2,x_3)  (1 - x_2) (1 -
  x_3) \, dx_2 \, dx_3 .
\end{split}
\end{equation}
We observe that $ J^{(1)}_{1,s_2,s_3} = - J^{(1)}_{0,s_2,s_3}$, and
that $J^{(2)}_{s_1,s_2,s_3}$ and $J^{(3)}_{s_1,s_2,s_3}$ can be
derived from $J^{(1)}_{s_1,s_2,s_3}$. So there are
three independent constants 
\begin{equation} \label{eq:define_J_S}
\begin{aligned}
  J_0 = {}& J_{0,0,0}
&
  J_1 = {}& J_{0,1,0}
&
  J_2 = {}& J_{0,1,1} .
\end{aligned}
\end{equation}

Due to the relations above, in the stiffness matrix each of the three
summand equals the tensor product of a 1-D discrete derivative (with
PML modifications), and a 2-D mass matrix (with PML modifiations).  We
first treat the PML modified derivative $- \pdpd{}{x_l} \alpha_l(x_l)
\pdpd{}{x_l}$.  Taking the case $l=1$, we can write the discrete
version of this as
\begin{equation}
  D^{(1)}_{2,{\rm FE},i,\tilde{i}}
  = \left\{ \begin{array}{ll}
 \frac{\alpha_{1,i-1/2}}{h_{1,i-1/2}} & \text{if $\tilde{i}=i-1$, }\\
 \frac{\alpha_{1,i+1/2}}{h_{1,i+1/2}}     & \text{if $\tilde{i}=i+1$, }\\
- \frac{\alpha_{1,i-1/2}}{h_{1,i-1/2}} - \frac{\alpha_{1,i+1/2}}{h_{1,i+1/2}}    
                                   & \text{if $\tilde{i}=i$, }\\
0 & \text{otherwise.}
\end{array}\right.
\end{equation}
The elements of the 2-D mass matrix with PML modifications read, for
the 2-D mass matrix related to the $(x_2,x_3)$ coordinate axes,
\begin{equation}
  N^{(2,3)}_{{\rm FE},j,k,\tilde{j},\tilde{k}}
  = J_{|\tilde{j}-j| + |\tilde{k}-k|}
  \sum_{(s_2,s_3) \in S(\tilde{j},j) \times S(\tilde{k},k)}
    \frac{h_{2,s_2} h_{3,s_3}}{\alpha_{2,s_2} \alpha_{3,s_3}} .
\end{equation}
when $\max( |\tilde{j}-j|, |\tilde{k}-k| ) \le 1$ (and is defined to
be 0 otherwise).  The full discrete Helmholtz operator becomes
\begin{equation} \label{eq:define_H_FE}
\begin{split}
  H_{{\rm FE},i,j,k,\tilde{i},\tilde{j},\tilde{k}}   
= {}& - M_{{\rm FE},i,j,k,\tilde{i},\tilde{j},\tilde{k}}
\\
  {}& - D^{(1)}_{2,{\rm FE},\tilde{i},i} N^{(2,3)}_{{\rm FE},j,k,\tilde{j},\tilde{k}}
  - D^{(2)}_{2,{\rm FE},\tilde{j},j} N^{(1,3)}_{{\rm FE},i,k,\tilde{i},\tilde{k}}
  - D^{(3)}_{2,{\rm FE},\tilde{k},k} N^{(1,2)}_{{\rm FE},i,j,\tilde{i},\tilde{j}} .
\end{split}
\end{equation}

\subsection{Coarse level optimized finite elements}
We will now show that the constants in the finite element method of
section~\ref{sec:discretization} can be chosen in such a way that the
rows associated with the regular, interior part of the mesh are equal
to the above described finite difference discretization, up to a
scalar factor. This means that the phase speeds of the coarse level
finite element method in the interior region closely match the phase
speeds of the fine level method. In this way we obtain the coarse
level discretization used in the two-grid method. The fine level
method is a finite difference method scaled by a constant 
$h^3$ (or $h^2$ in two dimensions), like in a finite element method.
We will start by assuming $k$ is constant.

Consider the expressions for the mass matrix $M_{{\rm FE}}$ and $N^{(l,m)}_{\rm FE}$. 
For the rows corresponding to degrees of freedom in the interior part
of the mesh, we have
\begin{equation}
  h_{1,i+s_1} = h_{2,j+s_2} = h_{3,k+s_3} = h 
\qquad
\text{and}
\qquad
  \alpha_{1,i+s_1} = \alpha_{2,j+s_2} = \alpha_{3,k+s_3} = 1 .
\end{equation}
since in the interior part of the mesh $\alpha_l = 1$ for $l=1,2,3$.
If we set 
\begin{equation} \label{eq:def_IJ_opt}
\begin{aligned}
  I_0 = {}& c_1/8 
&
  I_1 = {}& c_2/24
&
  I_2 = {}& c_3/24
\qquad\qquad
  I_3 = (1-c_1-c_2-c_3)/8
\\
  J_0 = {}& c_4/4
&
  J_1 = {}& c_5/8
&
  J_2 = {}& (1-c_4-c_5)/4 .
\end{aligned}
\end{equation}
then the operators defined in (\ref{eq:define_H_FE}) and
(\ref{eq:define_H_FDOPT}) have equal rows up to a factor $h^3$ in
three dimensions ($h^2$ in two dimensions).

The coarse scale finite element operator that we will consider is
given by taking (\ref{eq:def_IJ_opt}) as the definition of the $I_l$,
$l=0,1,2,3$ and $J_l$, $l=0,1,2$.

In appendix~\ref{app:FD_FE_appendix} we show that the shape function 
$\tilde{\psi}_{0,0,0}$ can be chosen such that the constants $I_l$ and 
$J_l$ satisfy the above equalities.

For variable $k$ we must specify how to obtain the coarse scale $k$
from the fine scale $k$. The coefficient $k$ at the coarse mesh cell
midpoints in PML layers are given by averaging with 
tensor products of 1-D averagings with 
$1/2, 1/2$ in the fine scale mesh points, and $1/4, 1/2, 1/4$ in the coarsened
interior part. The $\alpha$ values are evaluated at the cell-midpoints numerically.

For variable $k$ some differences between FD and FE discretizations
exist, due to the slightly different discretization of $k$ in these
operators.

\subsection{Finite element discretization with general test functions}
\label{app:FD_FE_appendix}

Equation (\ref{eq:def_IJ_opt}) contains a choice of the values $I_a$,
$a=0,1,2,3$ and $J_b$, $b=0,1,2$. Denote these prescribed values by
$\tilde{I}_a$ and $\tilde{J}_b$. We will show that there a shape
function $\tilde{\psi}_{0,0,0}$ such that the values of the $I_a$ and
$J_b$ defined in (\ref{eq:define_I_S}) and (\ref{eq:define_J_S}) 
agree with the prescribed values $\tilde{I}_a$ and $\tilde{J}_b$.

We define a symmetric 1-D tent function by
\begin{equation}
  T_{m,r}(x) = \left\{ 
    \begin{array}{ll} 
      (r - |x-m| )/r^2 & \text{if $|x-m|<r$} \\
      0                & \text{otherwise}
    \end{array} \right.
\end{equation}
for $m, r\in \R$, $r>0$. We define also define
\begin{equation}
  \tilde{T}_{r}(x) = \left\{ 
    \begin{array}{ll} 
      (r-x)/r   & \text{if $0 \le x \le r$} \\
      0         & \text{otherwise}
    \end{array} \right.
\end{equation}
Let $0< \eta$ be small, in each case $\eta < 1/2$, and let
$p_0 = \eta$, $p_1 = 1 - \eta$.
Given 7 parameters $A_a$, $B_b$, $a=0,1,2,3$ and $b=0,1,2$, we define 
$\tilde{\psi}_{0,0,0}$ by
\begin{equation} \label{eq:define_psitilde_000}
\begin{split}
  \tilde{\psi}_{0,0,0}
= {}&
    \sum_{i,j \in \{0,1\}} B_{i+j} T_{p_i,\eta}(x) T_{p_j,\eta}(y) \tilde{T}_{\eta}(z)
    + \sum_{i,k \in \{0,1\}} B_{i+k} T_{p_i,\eta}(x) T_{p_k,\eta}(z) \tilde{T}_{\eta}(y)
\\
{}&    + \sum_{j,k \in \{0,1\}} B_{j+k} T_{p_j,\eta}(y) T_{p_k,\eta}(z) \tilde{T}_{\eta}(x)
\\
{}& + \sum_{i,j,k \in \{0,1\}}
    \bigg( A_{i+j+k} - \frac{3-i-j-k}{2} B_{i+j+k} \eta \bigg)
    T_{p_i,\eta}(x) T_{p_j,\eta}(y) T_{p_k,\eta}(z)
\end{split}
\end{equation}

For the volume integrals $I_a$, we note that an approximation to $A_0
\delta$ is located at $(\eta,\eta,\eta)$, i.e.\ near $(0,0,0)$ and 
in the interior of the unit cube. Similarly,
approximate $\delta$ functions multiplied by one of the coefficients
$A_j$ are in all corners of the unit cube.
For the surface integrals $J_b$, we note that the restriction to the
plane $z=0$ contains an approximation to $B_0 \delta$ at
$(\eta,\eta)$ and similar approximations to $B_1 \delta$ and
$B_2 \delta$ in the other corners of the unit square.
The same is true for the planes $x=0$ and $y=0$.

Denote by $\Phi$ the linear map obtained by mapping
$  (A_0,A_1,A_2,A_3,B_0,B_1,B_2) $
to $\tilde{\psi}_{0,0,0}$
according to (\ref{eq:define_psitilde_000}) and then mapping 
$\tilde{\psi}_{0,0,0}$ to $(I_0,I_1,I_2,I_3,J_0,J_1,J_2)$ according to 
(\ref{eq:define_I_S}) and (\ref{eq:define_J_S}).

Let $\epsilon > 0$. We already observed that $\tilde{\psi}_{0,0,0}$ 
is a linear combination of approximate $\delta$
functions at the corners of the cube, supported just inside cube. 
This approximation becomes more accurate when $\eta \rightarrow 0$. 
Using this idea it is not difficult to show that when $\eta$ is
sufficiently small, then
\begin{equation}
\begin{split}
  | I_a - A_a | \le {}& \epsilon \| (A_0,A_1,A_2,A_3,B_0,B_1,B_2) \| 
, \qquad \text{for $a=0,1,2,3$}
\\
  | J_b - B_b | \le {}& \epsilon \| (B_0,B_1,B_2) \| 
, \qquad \text{for $b=0,1,2$}
\end{split}
\end{equation}
In other words, the linear map $\Phi$ is close to the identity,
we have $\| \Phi - I \| < C \epsilon$ (using the
matrix norm). This means that for sufficiently small $\eta$, the
linear map $\Phi$ is invertible and 
$  (A_0,A_1,A_2,A_3,B_0,B_1,B_2) $
can be found such that
\begin{equation}
(I_0,I_1,I_2,I_3,J_0,J_1,J_2)
=
(\tilde{I}_0,\tilde{I}_1,\tilde{I}_2,\tilde{I}_3,\tilde{J}_0,\tilde{J}_1,\tilde{J}_2) .
\end{equation}
Hence we have constructed $\tilde{\psi}_{0,0,0}$ with the desired property.

\end{document}